\documentclass[11pt,twoside]{amsart} 
\usepackage{amssymb}
\advance \topmargin by -\headheight
\advance \topmargin by -\headsep
\evensidemargin \oddsidemargin
\newtheorem{theorem}{Theorem}[section]
\newtheorem{corollary}[theorem]{Corollary}
\newtheorem{lemma}[theorem]{Lemma}
\newtheorem{proposition}[theorem]{Proposition}
\title{Semi-linear representations of $PGL$} 
\author{M.Rovinsky} 
\thanks{The author was supported in part by RFBR grant 02-01-22005}
\begin{document} 
\begin{abstract}
Let  $L$ be the function field of a projective space ${\mathbb P}^n_k$ 
over an algebraically closed field $k$ of characteristic zero, and $H$ 
be the group of projective transformations. An $H$-sheaf ${\mathcal V}$ 
on ${\mathbb P}^n_k$ is a collection of isomorphisms ${\mathcal V}
\longrightarrow g^{\ast}{\mathcal V}$ for each $g\in H$ 
satisfying the chain rule. 

We  construct, for any $n>1$, a fully faithful functor from the 
category of finite-dimensional $L$-semi-linear representations 
of $H$ extendable to the semi-group ${\rm End}(L/k)$ to the 
category of coherent $H$-sheaves on ${\mathbb P}^n_k$. 

The paper is motivated by a study of admissible representations 
of the automorphism group $G$ of an algebraically closed 
extension of $k$ of countable transcendence degree undertaken 
in \cite{rep}. The semi-group ${\rm End}(L/k)$ is considered 
as a subquotient of $G$, hence the condition on extendability. 

In the appendix it is shown that, if $\tilde{H}$ is either $H$, 
or a bigger subgroup in the Cremona group (generated by $H$ and 
a certain pair of involutions), then any semi-linear 
$\tilde{H}$-representation of degree one is an integral $L$-tensor power of 
$\det_L\Omega^1_{L/k}$. It is shown also  that this bigger subgroup 
has no non-trivial representations of finite degree if $n>1$. 
\end{abstract}
\maketitle 

\parskip 2mm 
\section{Introduction} 

Let $F$ be a field, $G$ a semigroup of endomorphisms of $F$ and $k=F^G$. 

An $F$-semi-linear $G$-representation is an $F$-space $V$ with a $k$-linear 
$G$-action such that $\sigma(a\cdot v)=\sigma a\cdot\sigma v$ for any 
$\sigma\in G$, $v\in V$ and $a\in F$. This is the same as a module over 
the associative central $k$-algebra $F\langle G\rangle:=
F\otimes_{{\mathbb Z}}{\mathbb Z}[G]$ with the evident 
left action of $F$ and the diagonal left action of $G$. We say 
that a semi-linear $G$-representation is {\sl non-degenerate} 
if the action of each element of $G$ is injective. 
%if $\dim_FV<\infty$ this is equivalent to $F\otimes_{\sigma(F)}\sigma(V)=V$ 
%If $G$ acts on $F$ faithfully, 
Semi-linear $G$-representations finite-dimensional 
over $F$ form an abelian tensor $k$-linear category. 
This category is rigid if the elements of $G$ are invertible. 
The set of isomorphism classes of non-degenerate semi-linear 
$F$-representations of $G$ of degree $r$ is canonically 
identified with the set $H^1(G,{\rm GL}_rF)$. 

\subsection{} \label{int-1} Let $k$ be a field of characteristic zero, 
and $k_{\infty}:=\cup_{j\ge 0}k_j$, where $k_0={\mathbb Q}$ and $k_{j+1}$ 
is generated over ${\mathbb Q}$ by all roots in $k$ of all elements 
of $k_j$. Suppose that $k$ contains all $\ell$-primary roots of unity 
for some prime $\ell\ge 2$. Let $L$ be the function field of a projective 
space ${\mathbb P}^n_k$ over $k$, and $G_n\cong{\rm PGL}_{n+1}k$ be the 
automorphism group of ${\mathbb P}^n_k$ over $k$. A $G_n$-{\sl structure} 
on a sheaf ${\mathcal V}$ on ${\mathbb P}^n_k$ is a collection of 
isomorphisms $\alpha_g:{\mathcal V}\stackrel{\sim}{\longrightarrow}
g^{\ast}{\mathcal V}$ for each $g\in G_n$ satisfying the chain rule: 
$\alpha_{hg}=g^{\ast}\alpha_h\circ\alpha_g$ for any $g,h\in G_n$. 
A $G_n$-{\sl sheaf} is a sheaf on ${\mathbb P}^n_k$ 
endowed with a $G_n$-structure. 

Let $\mathfrak{SL}_n$ be the full subcategory of the category of 
$L$-semi-linear representations of $G_n$, whose objects are restrictions 
of finite-dimensional $L$-semi-linear representations of the semi-group 
${\rm End}(L/k)$ to $G_n$. %\footnote{Let $Y$ be the complement of an 
%$(n+1)$-tuple of hyperplanes in ${\mathbb P}^n_k$ in general position, 
%and $H_n$ be the semi-subgroup in ${\rm End}(L/k)$ generated by its 
%subgroup $G_n$ and its sub-semisubgroup 
%${\rm End}_{{\rm dom}}(Y/k)$ of dominant endomorphisms of 
%$Y\cong({\mathbb G}_m)^n$. Clearly, $H_n$ is independent of $Y$. 
%In what follows we use only the sub-semigroup $H_n$ of ${\rm End}(L/k)$.}  

In this paper, for any $n\ge 2$, we construct a fully faithful functor 
$$\mathfrak{SL}_n\stackrel{{\mathcal S}}{\longrightarrow}
\{\mbox{coherent $G_n$-sheaves on ${\mathbb P}^n_k$}\}.$$ 
%and use it in Corollary \ref{podfakt} to show that any object of 
%$\mathfrak{SL}_n$ is a subquotient of the $L$-semi-linear 
%representation $W\otimes_kL$ of $G_n$ for a $k$-linear 
%representation $W$ of $G_n$ of finite degree. 

The construction of the functor ${\mathcal S}$ proceeds as follows. 
Fix a maximal split torus $T$ in $G_n$, and extend such a semi-linear 
representation $V$ to the semigroup ${\rm End}_{{\rm dom}}(Y/k)
\cong{\rm Mat}_{n\times n}^{\det\neq 0}{\mathbb Z}\ltimes T$ 
of dominant endomorphisms of the $n$-dimensional $T$-orbit 
$Y\cong({\mathbb G}_m)^n$ in ${\mathbb P}^n_k$. 
First we show (Proposition \ref{unip-tor}) that the restriction of $V$ to 
${\mathbb Z}_{\neq 0}^n\ltimes T$, where ${\mathbb Z}_{\neq 0}^n$ is a 
``maximal split torus'' in ${\rm Mat}_{n\times n}^{\det\neq 0}{\mathbb Z}$, 
is induced by a representation. An analytic argument (Lemma \ref{3-3}, 
we use here the assumption on existence of $\ell$-primary roots 
of unity in $k$) reduces this problem to a local result 
(Theorem \ref{main-loc}) asserting that any $k((t))$-semi-linear 
representation of the semi-group ${\mathbb N}$ (acting on the 
formal Laurent series field $k((t))$ by $p:t\mapsto t^p$) is 
induced by a representation. 

This implies (Lemma \ref{unip-rep}) that $V\longmapsto V^{T_{{\rm tors}}}$ 
gives a ``fibre functor'', i.e. $V=V^{T_{{\rm tors}}}\otimes_kL$, to 
the category of unipotent $k$-representations of $T$. Then, using more 
technical results of \S\ref{techno}, for each hyperplane $H$ stabilized 
by $T$ and for any $k_{\infty}$-lattice $U_0$ in the unipotent radical $U$ of 
the stabilizer $P$ of $H$, we construct (in Lemma \ref{restr-uni-rad}), 
another ``fibre functor'' $V\longmapsto V^{U_0}$ to the category of 
unipotent $k$-representations of $U$, so that the ${\mathcal O}
_{{\mathbb P}^n_k}({\mathbb P}^n_k-H)$-lattice ${\mathcal V}_H$ 
in $V$ spanned by $V^{U_0}$ is $P$-invariant and independent of $U_0$. 
%check that restriction of $V$ to $P$ is induced by 
%a representation trivial on (so trivial on $U$ if $k=k_{\infty}$). 

%Finally, (Lemma \ref{predp}) 
Localizing this lattice and varying $H$, one gets 
a coherent $G_n$-subsheaf ${\mathcal V}$ of the constant 
sheaf $V$ on ${\mathbb P}^n_k$ so that ${\mathcal V}|_Y=V^{T_{{\rm tors}}}
\otimes_k{\mathcal O}_Y=V^{U_0}\otimes_k{\mathcal O}_Y$ 
and $\Gamma({\mathbb P}^n_k-H,{\mathcal V})={\mathcal V}_H$. 

If $k=k_{\infty}$, one checks that the $G_n$-action on the total 
space $E$ of the vector bundle corresponding to ${\mathcal V}$ comes 
from a morphism of $k$-varieties $G_n\times_kE\longrightarrow E$, 
so the functor ${\mathcal S}$ factors through the 
category of $G_n$-equivariant coherent sheaves on ${\mathbb P}^n_k$, 
equivalent to the category of rational representations over $k$ 
of finite degree of the stabilizer of a point of ${\mathbb P}^n_k$. 

For $k$ transcendental over ${\mathbb Q}$ the objects of 
$\mathfrak{SL}_n$ are not equivariant sheaves anymore. 
%irreducible subquotients of are direct summands 
%of ${\rm Hom}_F((\Omega^n_{F/k})^{\otimes_F^i},
%(\Omega^1_{F/k})^{\otimes_F^j})$ for appropriate $i,j\ge 0$. 
For instance, there is a family of pairwise non-isomorphic 
semi-linear representations $\Omega^1_L/H\otimes_kL$ parametrized 
by the hyperplanes $H$ in the $k$-space $\Omega^1_k$. A choice of 
a non-zero element $v\in\Omega^1_k/H\cong k$ determines a non-split 
extension $0\longrightarrow L\stackrel{v\cdot}{\longrightarrow}
\Omega^1_L/H\otimes_kL\longrightarrow\Omega^1_{L/k}\longrightarrow 0.$ 
%for any $q\ge 1$, the Koszul resolution \begin{multline*}
%0\longrightarrow{\rm Sym}^q_k\Omega^1_k\otimes_kL\longrightarrow
%{\rm Sym}^{q-1}_k\Omega^1_k\otimes_k\Omega^1_L\longrightarrow\dots
%\longrightarrow
%{\rm Sym}^{q-p}_k\Omega^1_k\otimes_k\Omega^p_L\longrightarrow\\
%\dots\longrightarrow{\rm Sym}^2_k\Omega^1_k\otimes_k\Omega^{q-2}_L
%\longrightarrow\Omega^1_k\otimes_k\Omega^{q-1}_L\longrightarrow
%\Omega^q_L\longrightarrow\Omega^q_{L/k}\rightarrow 0 \end{multline*}
%gives rise to a natural map ${\rm Hom}_k({\rm Sym}^{q-p}_k\Omega^1_k,k)
%\stackrel{\alpha_{p,q}}{\longrightarrow}{\rm Ext}^{q-p}_{\mathfrak{SL}_m}
%(\Omega^q_{L/k},\Omega^p_{L/k})$, injective if $q=p+1=1$, where 
%$\mathfrak{SL}_m$ is an abelian category of $L$-semi-linear 
%representations of $G_n$ extendable to the semi-group ${\rm End}(L/k)$. 
%
%Clearly, ${\rm Hom}(\Omega^i_L,\Omega^j_{L/k})$ is 
%one-dimensional over $k$ if $i=j$, and it is zero otherwise.  
%$${\rm Hom}(\Omega^i_L,\Omega^j_L)=\left\{\begin{array}{ll} 
%\Omega^{j-i}_k & \mbox{if $j\ge i$} \\ 0 & \mbox{otherwise} 
%\end{array}\right.$$

\subsection{} The paper is motivated by a study of admissible 
representations\footnote{A representation of a topological group 
is called {\sl smooth} if the stabilizers are open. A smooth 
representation of is called {\sl admissible} if the subspace 
fixed by any open subgroup is finite-dimensional.} of 
the automorphism group $G$ of an algebraically closed 
extension $F$ of $k$ of countable transcendence degree 
undertaken in \cite{rep}. One can expect that any such representation 
is contained in appropriate admissible semi-linear representation, 
cf. \S\ref{conj-cor}. In \S\ref{pgl-systems} an abelian category 
${\mathcal P}$ of compatible systems of semi-linear representations 
of Cremona groups $Cr_n(k)$ is introduced and a faithful functor %$PGL_n$'s
$\{\mbox{smooth $F$-semi-linear representations of $G$}\}
\stackrel{\Phi}{\longrightarrow}{\mathcal P}$ is 
constructed. 

It is well-known that any irreducible $G_n$-equivariant 
coherent sheaf on ${\mathbb P}^n_k$ is a direct summand of 
${\rm Hom}_{{\mathcal O}_{{\mathbb P}^n_k}}
((\Omega^n_{{\mathbb P}^n_k/k})^{\otimes_{{\mathcal O}_{{\mathbb P}^n_k}}^i},
\bigotimes_{{\mathcal O}_{{\mathbb P}^n_k}}^{\bullet}
\Omega^1_{{\mathbb P}^n_k/k})$ for an appropriate $i\ge 0$. 

Then Theorem \ref{main-th} and the results presented in 
\S\ref{main-conj} suggest that irreducible admissible 
representations of $G$ are contained in the algebra of 
relative differential forms $\Omega^{\bullet}_{F/k}$. 
If the same is true for any irreducible object of ${\mathcal I}_G$ 
(cf. \S\ref{conj-cor}, p.\pageref{conj-cor}) then for any smooth 
proper variety with no regular differential forms of degree 
$\ge 2$ the Albanese map identifies the group of classes of 
0-cycles of degree zero modulo rational equivalence with 
the Albanese variety, cf. Corollary \ref{bl-con}. 
This would also imply a description of pure motives as 
admissible $G$-modules, cf. Corollary \ref{adm-impl-adm}. 

\subsection{} In the appendix semi-linear representations of degree 
one are studied in more detail. The main results there are Corollary 
\ref{pgl-2} and Proposition \ref{cre}, where it is shown that if $L$ 
is the function field of a projective space ${\mathbb P}^n_k$ over 
an algebraically closed field $k$ of characteristic zero and $G$ is 
either the group of projective transformations, or a certain bigger 
subgroup in the Cremona group, then any semi-linear $G$-representation 
of degree one is an integral $L$-tensor power of $\det_L\Omega^1_{L/k}$. 
This bigger subgroup has such an advantage that it has no non-trivial 
representations of finite degree if $n\ge 2$ (cf. Proposition 
\ref{no-rep}), so at least this source of unexpected 
semi-linear representations is excluded. 

\section{Examples of semi-linear $G$-representations} \label{simpl-ex}
In the first example $G$ is a (semi-)group generated by a single 
element $T$. Such situation (in a greater generality) was studied 
in \cite{O}. Let $V$ be an $F$-semi-linear 
$G$-representation of degree $N$ admitting a cyclic vector $v$, 
i.e., such that $\{v,Tv,\dots,T^{N-1}v\}$ is an $F$-base of $V$. 
Then $T^Nv=\sum_{j=0}^{N-1}h_jT^jv$ for some $h_0,\dots,h_{N-1}\in F$, 
so the matrix of $T$ in this base looks as $\left(\begin{array}{cccccc}
0 & 0 & \dots & 0 & 0 & h_0 \\ 1 & 0 & \dots & 0 & 0 & h_1 \\ 
\dots & \dots & \dots & \dots & \dots & \dots \\ 
0 & 0 & \dots & 0 & 1 & h_{N-1} \end{array}\right)$. 

The following well-known result shows that this situation is typical. 
\begin{lemma} \label{cyclicity} 
Let $\sigma$ be an endomorphism of a field $F$ of infinite order. 
Then any non-degenerate $F$-semi-linear representation of $\sigma$ 
of finite degree is cyclic. \end{lemma}
{\it Proof.} We proceed by induction on degree (or length) $N$, the 
case $N=1$ being trivial. For $N>1$ let $0\neq V_0\subset V$ be an 
irreducible subrepresentation. By induction assumption, $V/V_0$ is 
generated by some $\overline{v}\in V/V_0$. Choose a lift $v\in V$ 
of $\overline{v}$. 

Suppose that $V$ is not cyclic. Then $V=V_0\oplus
\langle v+w\rangle$ for any $w\in V_0$. Then the left ideal $Ann(v)$ 
is contained in the two-sided ideal $Ann(V_0)$. If a non-zero 
two-sided ideal contains $\sigma^m+f_{m-1}\sigma^{m-1}+\dots+f_0$ 
with minimal possible $m\ge 0$ then it contains also $\sigma^m+
f_{m-1}\frac{\sigma^{m-1}\lambda}{\sigma^m\lambda}\sigma^{m-1}+
\dots+f_0\frac{\lambda}{\sigma^m\lambda}$ for any $\lambda\in 
F^{\times}$, so $f_{m-1}=\dots=f_0=0$. As $\dim_F\langle v\rangle
=\dim_F(F\langle\sigma\rangle/Ann(v))\le N<\infty$, one has 
$Ann(v)\neq 0$. In particular, $Ann(V_0)\neq 0$ and $\sigma^m(V_0)=0$ 
for some $m\ge 0$, which contradicts our assumptions, and thus, $V$ 
is cyclic. \qed 

\subsection{Linear and semi-linear representations} 
\label{flat-source}
To any subfield $k'$ in $F$ invariant under the $G$-action and 
to any finite-dimensional semi-linear $k'$-representation 
$V_0$ of $G$ one associates the semi-linear $F$-representation 
$V_0\otimes_{k'}F$ of $G$. On the level of isomorphism 
classes of semi-linear representations of $G$ of degree 
$r$ this operation coincides with the natural map 
\begin{equation} \label{ext-coeff} H^1(G,{\rm GL}_rk')
\longrightarrow H^1(G,{\rm GL}_rF).\end{equation} 

The following lemma gives a sufficient condition 
for injectivity of the map (\ref{ext-coeff}). 
\begin{lemma} \label{source} Let $k'$ be a Galois extension of $k$ in $F$. 
If the $G$-orbit of any element of $F-k'$ spans a $k'$-subspace in $F$ of 
dimension $>r^2$ then the map (\ref{ext-coeff}) is injective. \end{lemma}
{\it Proof.} Let $(a_{\sigma})$ and $(a'_{\sigma})$ be two 1-cocycles 
representing some classes in $H^1(G,{\rm GL}_rk')$. Suppose that they 
become the same in $H^1(G,{\rm GL}_rF)$, i.e., there is an element 
$b\in{\rm GL}_rF$ such that $a_{\sigma}=b^{-1}a'_{\sigma}\sigma b$ 
for all $\sigma\in G$. Equivalently, $ba_{\sigma}=a'_{\sigma}\sigma b$ 
for all $\sigma\in G$. If $b\not\in{\rm GL}_rk'$, i.e., there are some 
$1\leq s,t\leq r$ such that $b_{st}\not\in k'$, then there is 
$\sigma\in G$ such that 
$\sigma b_{st}\not\in\langle b_{ij}~|~1\leq i,j\leq r\rangle_{k'}$, 
which contradicts $ba_{\sigma}=a'_{\sigma}\sigma b$. This means 
that $b\in{\rm GL}_rk'$, and thus, the classes of $(a_{\sigma})$ 
and $(a'_{\sigma})$ in $H^1(G,{\rm GL}_rk')$ coincide. \qed 

\vspace{5mm}

In opposite direction, let $V$ be an $F$-semi-linear 
$G$-representation, and $\rho:G\longrightarrow{\rm GL}(V)$ 
an $F$-linear representation. Set 
$V_{\rho}:=\{w\in V~|~\sigma w=\rho(\sigma)w\}$. 
Then $V_{\rho}$ is a $k$-space. 
\begin{lemma} The natural map 
$F\otimes_kV_{\rho}\longrightarrow V$ is injective. \end{lemma}
{\it Proof.} Let $\{w_1,\dots,w_m\}\subset V_{\rho}$ be linear 
independent over $F$. Suppose that $w=\sum_j\lambda_jw_j\in V_{\rho}$. 
Then $\sigma w-\rho(\sigma)w=\sum_j(\sigma\lambda_j-\lambda_j)
\rho(\sigma)w_j=0$, and therefore, $\sigma\lambda_j=\lambda_j$ 
for any $j$, i.e., $\lambda_j\in k$. \qed 

\vspace{5mm}

{\it Remark.} The irreducibility of a representation $W$ of $G$ 
over $k$ does not imply the irreducibility of the $F$-semi-linear 
$G$-representation $F\otimes_kW$. 

For example, let $Q$ be a finite-dimensional $k$-space, 
$F=k({\mathbb P}(Q))$ be the function field of its projectivization, 
$G={\rm PGL}(Q)$, and $W=\mathfrak{sl}(Q)$ be the adjoint 
representation of $G$. We identify $W$ with the global vector fields 
on ${\mathbb P}(Q)$, thus getting a non-injective surjection 
$W\otimes_kL\longrightarrow{\rm Der}(L/k)$. 

\section{Some semi-linear representations 
of groups exhausted by finite subgroups}
Let $F$ be a field, $G$ be a group of field automorphisms of 
$F$ and $H_1\triangleleft H_2\subseteq G$. Set $L:=F^{H_1}$. 
Suppose that $H:=H_2/H_1$ is exhausted by its finite subgroups, 
i.e., $H$ is a torsion group and any pair of its 
finite subgroups generates a finite subgroup. 

By Hilbert Theorem 90, one has $$Z^1(H,GL_NL)=
{\lim_{\longleftarrow}}_{\Phi}GL_NL^{\Phi}\backslash GL_NL,
\quad\mbox{and}\quad H^1(H,GL_NL)=Z^1(H,GL_NL)/GL_NL,$$
where $\Phi$ runs over the set of finite subgroups in $H$. 

\subsection{Endomorphisms and contractions} \label{contra}
Suppose that $\xi\in G$ induces an endomorphism of $H$, i.e., 
$\xi^{-1}H_2\xi\subseteq H_2$ and $\xi^{-1}H_1\xi\subseteq H_1$. 
Given an $F$-semi-linear $G$-representation $V$, consider $V^{H_1}$. 
This is an $L$-semi-linear $H$-representation with an injective 
semi-linear $\xi$-action: for any $h\in H_1$ and any $v\in V^{H_1}$ 
one has $h(\xi v)=\xi(\xi^{-1}h\xi)v=\xi v$, so $\xi v\in V^{H_1}$. 

Let $N=\dim_LV^{H_1}$. For any $\zeta\in H$ one has $f_{\zeta}\cdot
\zeta f_{\xi}=f_{\xi}\cdot\xi f_{\xi^{-1}\zeta\xi}\in{\rm GL}_NL$. 
As $f_{\zeta}=f^{-1}\cdot\zeta f$ for some $f\in
{\lim\limits_{\longleftarrow}}_{\Phi}{\rm GL}_NL^{\Phi}
\backslash{\rm GL}_NL$, this
implies that $(\xi f)^{-1}\cdot\zeta(\xi f)=(f\cdot f_{\xi})^{-1}
\cdot\zeta(f\cdot f_{\xi})$, so $f\cdot f_{\xi}=\xi f\in
{\lim\limits_{\longleftarrow}}_{\Phi}{\rm GL}_NL^{\Phi}\backslash{\rm GL}_NL$. 

{\sc Example.} Suppose that $\xi$ contracts $H$, i.e., for any 
$\zeta\in H_2$ there is $s\ge 1$ such that $\xi^{-s}\zeta\xi^s\in H_1$. 
Then $\zeta(\xi^sg)=\xi^s(\xi^{-s}\zeta\xi^s)g=\xi^sg$, i.e.,  
$\xi^sg\in{\rm GL}_NL^{\langle\zeta\rangle}$ for any $g\in{\rm GL}_NL$. 
This implies that for any $f\in Z^1(H,{\rm GL}_NL)$ as above there 
exists the limit $\lim\limits_{s\to\infty}\xi^sf\cdot f_{\xi^s}^{-1}=
\lim\limits_{s\to\infty}f_{\xi^s}^{-1}=\cdots\xi^3 f_{\xi}^{-1}\cdot
\xi^2 f_{\xi}^{-1}\cdot\xi f_{\xi}^{-1}\cdot f_{\xi}^{-1}=f$, so 
there is a bijection ${\rm GL}_NL\stackrel{\sim}{\longrightarrow}Z^1
(H\rtimes\langle\xi\rangle,{\rm GL}_NL)$ sending $g\in{\rm GL}_NL$ 
to $f_{\xi}=g$, $f_{\zeta}=f_{\xi^s}\cdot\zeta f_{\xi^s}^{-1}=
f^{-1}\cdot\zeta f$ if $\zeta\in H$ and $\xi^s(\zeta)=1$. 

%In view of this example, it is natural to consider non-degenerate 
%semi-linear representations of semigroups of endomorphisms of fields. 

\section{Local problem: $k((t))$-semi-linear ${\mathbb N}$-representations}
Let ${\mathbb N}$ be the multiplicative semigroup of positive 
integers acting on $k((t))$ by $p:t\mapsto t^p$. In this section 
we show that any semi-linear representation of ${\mathbb N}$ 
finite-dimensional over $k((t))$ is induced by a 
$k$-representation. 
\begin{lemma} \label{integ} % One has 
$H^1(S,{\rm GL}_Nk)\stackrel{\sim}{\longrightarrow}
H^1(S,{\rm GL}_Nk[[t]])\hookrightarrow H^1(S,{\rm GL}_Nk((t)))$ 
and $H^1(S,1+t\mathfrak{gl}_Nk[[t]])=\{\ast\}$ 
for any subsemigroup $S\subseteq{\mathbb N}$. \end{lemma}
{\it Proof.} We may suppose that $S\neq\{1\}$. For some 
$p\in S-\{1\}$ and $(f_{\ell})\in Z^1(S,{\rm GL}_Nk[[t]])$. 
set $B:=f_p^{-1}-f_p(0)^{-1}\in t\mathfrak{gl}_Nk[[t]]$. 
As $f_p(0)^{s+1}f_{p^{s+1}}^{-1}=f_p(0)^{s+1}(f_p(0)^{-1}+"{p^s}"B)
f_{p^s}^{-1}=f_p(0)^sf_{p^s}^{-1}+f_p(0)^{s+1}\cdot p^sB\cdot 
f_{p^s}^{-1}$, there is the limit $\Phi:=\lim\limits_{s\to\infty}
f_{p^s}(0)f_{p^s}(t)^{-1}=\lim\limits_{s\to\infty}f_p(0)^sf_{p^s}^{-1}
\in{\bf 1}_N+t\mathfrak{gl}_Nk[[t]]$. Then 
$\Phi(t)f_p(t)\Phi(t^p)^{-1}=f_p(0)$, so we may suppose that 
$f_p\in{\rm GL}_Nk$ and $f_p={\bf 1}_N$ if $f_p\equiv{\bf 1}_N\pmod t$. 
As $f_p^{-1}f_{\ell}(t)f_p=f_{\ell}(t^p)$, this implies that 
$f_{\ell}\in{\rm GL}_Nk$ for any $\ell\in S$ and 
$f_{\ell}={\bf 1}_N$ if $f_{\ell}\equiv{\bf 1}_N\pmod t$. 
The injectivity statement follows from a similar argument. \qed 

\begin{lemma} \label{one-dim} ${\rm Hom}(S,k^{\times})\times 
d(S)^{-1}{\mathbb Z}/{\mathbb Z}\stackrel{\sim}{\longrightarrow}
H^1(S,k((t))^{\times})$ for any non-trivial subsemigroup 
$S\subseteq{\mathbb N}$, where $d(S)$ is the greatest 
common divisor of $s-1$ for all $s\in S$. \end{lemma}
{\it Proof.} By Lemma \ref{integ}, one has $H^1(S,1+tk[[t]])=0$. 
Clearly, $H^1(\langle p\rangle,t^{{\mathbb Z}})\cong
{\mathbb Z}/(p-1){\mathbb Z}$ for any $p\ge 2$. 
As $k((t))^{\times}\cong{\mathbb Z}\times k^{\times}
\times(1+tk[[t]])$, one has $H^1(\langle p\rangle,k((t))^{\times})
\cong{\mathbb Z}/(p-1){\mathbb Z}\times k^{\times}$. 
Then we may assume that $f_p\in t^{m_p}k^{\times}$ for some 
$m_p\in{\mathbb Z}$, so $f_{\ell}(t^p)=f_{\ell}(t)t^{(\ell-1)m_p}$. 
As $f_{\ell}=a_{\ell}t^{m_{\ell}}\varphi_{\ell}$ for some $a_{\ell}
\in k^{\times}$, $m_{\ell}\in{\mathbb Z}$ and $\varphi_{\ell}\in 
1+tk[[t]]$, we conclude that $\varphi_{\ell}=1$ and 
$\frac{m_{\ell}}{\ell-1}=\frac{m_p}{p-1}$, i.e., that 
$\frac{m_{\ell}}{\ell-1}\in\frac{1}{d(S)}{\mathbb Z}$ 
does not depend on $\ell$. \qed 

\begin{lemma} \label{block-tri} If $f(t)\in\mathfrak{gl}_Nk[[t]]
\cap{\rm GL}_Nk((t))$ and $\ell\ge 2$ then there is an element 
$g(t)\in{\rm GL}_Nk[[t]]$ such that the matrix $g(t)^{-1}f(t)g(t^{\ell})$ 
is blockwise upper triangular with at most two diagonal blocks: one 
invertible and constant, and another nilpotent modulo $t$. \end{lemma}
{\it Proof.} Fix an element $A\in{\rm GL}_Nk$ sending first 
${\rm rk}f(0)^N$ coordinate vectors to ${\rm Im}f(0)^N$ and the rest 
to $\ker f(0)^N$. Then the matrix $A^{-1}f(0)A$ is blockwise diagonal 
with two non-zero blocks, the first one invertible and the second one 
nilpotent. So we may suppose that $f(t)=\left(\begin{array}{cc} 
E & F \\ G & H \end{array}\right)$, where $E\in{\rm GL}_mk[[t]]$, 
$0\le m\le N$, $H\in\mathfrak{gl}_{N-m}k[[t]]$ is nilpotent modulo 
$t$, $F\in t{\rm Mat}_{m\times(N-m)}k[[t]]$ and 
$G\in t{\rm Mat}_{(N-m)\times m}k[[t]]$. 

Let $C_0(t)=G(t)E(t)^{-1}\in t{\rm Mat}_{(N-m)\times m}k[[t]]$ and 
$C_{j+1}(t)=(G(t)+H(t)C_j(t^{\ell}))(E(t)+F(t)C_j(t^{\ell}))^{-1}$ 
for any $j\ge 0$. By induction on $j\ge 1$ we check that $C_j
\equiv C_{j-1}\pmod{t^{\ell^j}}$. 
For $j=1$ this follows from $t^{\ell}|C_0(t^{\ell})$. 

Suppose now that $C_j\equiv C_{j-1}\pmod{t^{\ell^j}}$ for some $j\ge 1$. 
Then \begin{multline*}C_j(t^{\ell})\equiv C_{j-1}(t^{\ell})\pmod{
t^{\ell^{j+1}}},\quad\mbox{and therefore,}\\ 
(G(t)+H(t)C_j(t^{\ell}))(E(t)+F(t)C_j(t^{\ell}))^{-1}\\ 
\equiv(G(t)+H(t)C_{j-1}(t^{\ell}))(E(t)+F(t)C_{j-1}(t^{\ell}))^{-1}
\pmod{t^{\ell^{j+1}}},\end{multline*} 
which is equivalent to $C_{j+1}\equiv C_j\pmod{t^{\ell^{j+1}}}$. 
This implies that the sequence $(C_j)_{j\ge 0}$ is convergent in 
$t{\rm Mat}_{(N-m)\times m}k[[t]]$. Denote by $C(t)$ its limit. 
Then $C(t)F(t)C(t^{\ell})=G(t)-C(t)E(t)+H(t)C(t^{\ell})$, and therefore, 
\begin{multline*} \left(\begin{array}{cc} E(t)+F(t)C(t^{\ell}) & F(t) \\ 
G(t)-C(t)E(t)+H(t)C(t^{\ell})-C(t)F(t)C(t^{\ell}) & H(t)-C(t)F(t) 
\end{array}\right) \\ 
=\left(\begin{array}{cc} E'(t) & F(t) \\ 0 & H'(t) \end{array}\right)
=\left(\begin{array}{cc} {\bf 1}_m & 0 \\ -C(t) & {\bf 1}_{N-m} 
\end{array}\right)f(t)\left(\begin{array}{cc} {\bf 1}_m & 0 \\ 
C(t^{\ell}) & {\bf 1}_{N-m} \end{array}\right) 
\end{multline*} 
where $E'\in{\rm GL}_mk[[t]]$ and $H'$ 
is nilpotent modulo $t$, since $C\equiv 0\pmod t$. 
The rest follows from Lemma \ref{integ}. \qed

\begin{lemma} \label{diag-to-cont} Let $S\subseteq{\mathbb N}$ be 
a subsemigroup containing a pair of coprime elements $q,\ell\ge 2$. 
If $f_p(t)\in{\rm GL}_Nk((t))$ are upper triangular and 
$f_p(t)f_{\ell}(t^p)=f_{\ell}(t)f_p(t^{\ell})$ for all $p\in S$, 
then there is an element $g(t)\in{\rm GL}_Nk((t))$ such that 
$g(t)^{-1}f_p(t)g(t^p)\in{\rm GL}_Nk$ for all $p\in S$. \end{lemma} 
{\it Proof.} We proceed by induction on $N\ge 1$, the case $N=1$ 
follows from Lemma \ref{one-dim}. Let $f_p(t)=\left(\begin{array}{cc} 
A_p(t) & B_p(t) \\ 0 & D_p(t) \end{array}\right)$. By induction 
assumption, acting by a blockwise diagonal coboundary, we may
suppose that $A_p\in{\rm GL}_{N-1}k$ and $D_p\in k^{\times}$. The 
commutativity condition becomes $A_pB_{\ell}(t^p)+B_p(t)D_{\ell}=
A_{\ell}B_p(t^{\ell})+B_{\ell}(t)D_p$, or $c_p\varphi_{\ell}(t^p)
+\varphi_p(t)=c_{\ell}\varphi_p(t^{\ell})+\varphi_{\ell}(t)$, 
where $c_p=A_pD_p^{-1}$ and $\varphi_p(t)=B_p(t)D_p^{-1}$. 
If we replace $f_p(t)$ by $f'_p(t)=\left(\begin{array}{cc} 
{\bf 1}_{N-1} & -B \\ 0 & 1 \end{array}\right)f_p\left(\begin{array}{cc} 
{\bf 1}_{N-1} & B(t^p) \\ 0 & 1 \end{array}\right)=\left(
\begin{array}{cc} A_p & A_pB(t^p)+B_p-BD_p \\ 0 & D_p 
\end{array}\right)$ then $\varphi'_{\ell}=\varphi_{\ell}+c_{\ell}
B(t^{\ell})-B$. After an appropriate choice of $B\in k[t^{-1}]$ 
(as $ht^{-\ell m}$ is `cohomologeous to $c^{-1}_{\ell}ht^{-m}$'), 
we may assume that $\varphi_{\ell}\in
(\sum_{j=1}^{\ell-1}t^{-j}k[t^{-{\ell}}]+k[[t]])^{N-1}$. 

Let $t^m\varphi_{\ell}\in k[[t]]^{N-1}$ for minimal possible $m$. 
Suppose that $m\ge 1$. Then $m$ is not divisible by $\ell$. Let 
$p\in S-\{1\}$ be prime to $\ell$. By induction on $s\ge 0$, one checks 
that $\varphi_p(t)=\Phi_s(t)+c_{\ell}^{s+1}\varphi_p(t^{\ell^{s+1}})$, 
where $\Phi_s(t)=\sum_{j=0}^sc_{\ell}^j(\varphi_{\ell}(t^{\ell^j})-
c_p\varphi_{\ell}(t^{\ell^jp}))$. Clearly, $k[[t]]^{N-1}$ contains 
$t^{\ell^spm}\Phi_s(t)$, but does not contain $t^{\ell^spm-1}\Phi_s(t)$. 
As $\ell^{s+1}$ does not divide $\ell^spm$, one has also 
$t^{\ell^spm-1}\varphi_p(t)\not\in k[[t]]^{N-1}$, 
which is impossible for $s\gg 0$.

This contradiction shows that $\varphi_{\ell}\in k[[t]]^{N-1}$, 
i.e., that $f_{\ell}\in{\rm GL}_Nk[[t]]$. 

By Lemma \ref{integ}, we may suppose that $f_{\ell}\in{\rm GL}_Nk$. 
Then the commutativity condition becomes $f_{\ell}f_p(t^{\ell})=
f_p(t)f_{\ell}$, and thus, $f_p(t)\in{\rm GL}_Nk$ for all $p\in S$. \qed

\begin{theorem} \label{main-loc} Let $N\ge 1$ be an integer, 
$S\subseteq{\mathbb N}$ be a subsemigroup containing a pair of 
integers $\ell\ge p\ge 2$ such that the least common multiple 
$[p-1,\dots,p^N-1]$ divides $\ell$, and containing a pair 
of coprime integers $\ge 2$. Then $H^1(S,{\rm GL}_Nk)
\stackrel{\sim}{\longrightarrow}H^1(S,{\rm GL}_Nk((t)))$. \end{theorem} 
{\it Proof.} We proceed by induction on $N\ge 1$, the case $N=1$ 
being contained in Lemma \ref{one-dim}. For $N>1$ let $V$ be a 
semi-linear representation of $S$, which is $N$-dimensional over 
$F:=k((t))$. By Lemma \ref{cyclicity}, there is a vector $v\in V$ 
generating $V$ as $F\langle\sigma\rangle$-module, where 
$\sigma:=\sigma_p$. Let $\sigma^Nv=\sum_{j=0}^{N-1}h_j\cdot\sigma^jv$, 
where $h_j\in F$, and $\alpha=\max_{0\le j<N}\frac{v(h_j)}{p^j-p^N}\in
\frac{1}{[p-1,\dots,p^N-1]}{\mathbb Z}$. Set $v'=t^{p^{N-1}\ell\alpha}
\cdot\sigma^{N-1}\tau v$, where $\tau=\sigma_{\ell}$. Then 
$$\sigma^Nv'%=t^{p^{2N-1}\ell\alpha}\sigma^{N-1}\tau(\sigma^Nv)
=t^{p^{2N-1}\ell\alpha}\sum_{j=0}^{N-1}\sigma^{N-1}\tau h_j
\cdot\sigma^j(\sigma^{N-1}\tau v)=\sum_{j=0}^{N-1}
t^{p^{N-1}(p^N-p^j)\ell\alpha}\sigma^{N-1}\tau h_j\cdot\sigma^jv'.$$

Set $h'_j=t^{p^{N-1}(p^N-p^j)\ell\alpha}\sigma^{N-1}\tau h_j$. 
Then $v(h_j')=p^{N-1}(p^N-p^j)\ell(\alpha+\frac{v(h_j)}{p^N-p^j})\ge 0$, 
and $\min_{0\le j<N}v(h_j')=0$. 
This means, that in the basis $\{v'\sigma v',\dots,\sigma^{N-1}v'\}$ 
the matrix of $\sigma$ is $$\left(\begin{array}{cccccc}
0 & 0 & \dots & 0 & 0 & h_0' \\ 1 & 0 & \dots & 0 & 0 & h_1' \\ 
\dots & \dots & \dots & \dots & \dots & \dots \\ 
0 & 0 & \dots & 0 & 1 & h_{N-1}' \end{array}\right),$$ 
where $h_j'\in k[[t]]$ and $h_m'\in k[[t]]^{\times}$ for some 
$0\le m<N$, i.e., the matrix of $\sigma$ is not nilpotent at 0. 
Let $m$ be minimal. According to Lemma \ref{block-tri}, in an 
appropriate basis, the matrix of $\sigma$ belongs to 
$\left(\begin{array}{cc} {\rm GL}_{N-m}k & {\rm Mat}k[[t]] \\
0 & \mathfrak{gl}_mk[[t]] \end{array}\right)$, i.e, 
there exists a non-zero $k$-subspace $W_0\subset V$ invariant under 
$\sigma$ and such that the natural map $W_0\otimes_kF\longrightarrow V$ 
is injective, and therefore, there is a non-constant polynomial 
$P\in k[T]$ such that $\ker P(\sigma)\neq 0$. Choose such $P$ 
with minimal possible degree (in particular, of degree one if 
$k$ is algebraically closed). Then the natural map 
$\ker P(\sigma)\otimes_kF\longrightarrow V$ is injective. If 
$P(\sigma)v=0$ then $P(\sigma)\xi v=0$ for any $\xi\in S$, so 
$\ker P(\sigma)$ is $S$-invariant. This implies that there is a 
basis, where the matrices of all elements of $S$ belong to 
$\left(\begin{array}{cc} {\rm GL}_sk & {\rm Mat}_{s\times(N-s)}k((t)) \\
0 & {\rm GL}_{N-s}k((t)) \end{array}\right)$, where 
$1\le s=\dim\ker P(\sigma)\le N$. By induction hypothesis, applied 
to $V/\ker P(\sigma)\otimes_kF$, there is a basis, 
where the matrices of all elements of $S$ belong to 
$\left(\begin{array}{cc} {\rm GL}_sk & {\rm Mat}_{s\times(N-s)}k((t)) \\
0 & {\rm GL}_{N-s}k \end{array}\right)$. Over $\overline{k}$, 
the matrices of all elements can be made upper triangular. 
It follows from Lemma \ref{diag-to-cont} that over $\overline{k}$ 
the matrices of all elements can be made constant, and therefore, 
the same can be done over $k$. \qed

\section{Purely transcendental extensions: reduction to the local 
problem} 
\begin{lemma} \label{red-to-1} For each $1\le j\le n$ fix a 
sub-semigroup $H_j$ in ${\rm End}(k(t)/k)$ containing $t\mapsto 
t^{\ell_j}$ for some $\ell_j\ge 2$. Suppose that the natural map 
${\rm Hom}(H_j,{\rm GL}_N\overline{K})\longrightarrow H^1
(H_j,{\rm GL}_N\overline{K}(t))$ is surjective for any extension 
$K$ of $k$, where ${\rm End}(k(t)/k)$ is considered as sub-semigroup 
of ${\rm End}(K(t)/K)$. Then ${\rm Hom}(\prod_{j=1}^nH_j,{\rm GL}_NK)
\longrightarrow H^1(\prod_{j=1}^nH_j,{\rm GL}_NK(x_1,\dots,x_n))$ 
is surjective for any $n\ge 1$. \end{lemma} 
{\it Proof.} Assume first that $n=1$. Let $(f_{\sigma})\in 
Z^1(H,{\rm GL}_NK)$. Then $f_{\sigma}=f_1(x)^{-1}
\cdot g_{\sigma}\cdot f_1(\sigma x)$, where $g_{\sigma}\in
{\rm GL}_N\overline{K}$ and $f_1(x)\in{\rm GL}_N\overline{K}(x)$ 
for any $\sigma\in H$. Fix some $\alpha\in k$ such that $f_1(\alpha)$ 
is a well-defined element of ${\rm GL}_N\overline{K}$. Set $f(x):=
f_1(\alpha)^{-1}f_1(x)$ and $h_{\sigma}:=f_1(\alpha)^{-1}
g_{\sigma}f_1(\alpha)$ (so $h_{\sigma}\in{\rm GL}_N\overline{K}$). 
Then $f_{\sigma}(x)=f(x)^{-1}\cdot h_{\sigma}\cdot f(\sigma x)$ and 
$f(\alpha)=1$. As $f_{\sigma}\in{\rm GL}_NK(x)$, we get 
$f^{\tau}(x)^{-1}\cdot h_{\sigma}^{\tau}\cdot f^{\tau}(\sigma x)=
f(x)^{-1}\cdot h_{\sigma}\cdot f(\sigma x)$ for any $\tau\in
{\rm Gal}(\overline{K}(x)/K(x))={\rm Gal}(\overline{K}/K)$, which 
is equivalent to $f(x)\cdot f^{\tau}(x)^{-1}=h_{\sigma}\cdot
(f(\sigma x)f^{\tau}(\sigma x)^{-1})(h_{\sigma}^{\tau})^{-1}$. 
Looking at $\sigma$'s of type $x\mapsto x^{\ell^q}$  we see that 
$f(x)\cdot f^{\tau}(x)^{-1}\in\bigcap_{q\ge 1}{\rm GL}_N\overline{K}
(x^{\ell^q})={\rm GL}_N\overline{K}$. As $f(\alpha)=1$, one has 
$f^{\tau}(\alpha)=1$, and thus, $f(x)=f^{\tau}(x)$ for any 
$\tau\in{\rm Gal}(\overline{K}/K)$, so $f(x)\in{\rm GL}_NK(x)$, 
and therefore, $h_{\sigma}\in{\rm GL}_NK(x)\cap{\rm GL}_N\overline{K}
={\rm GL}_NK$. This means that the natural map 
${\rm Hom}(H,{\rm GL}_NK)\longrightarrow H^1
(H_j,{\rm GL}_NK(x))$ is surjective. 

Suppose now that $n\ge 2$. We proceed by induction on $n$. Consider 
a cocycle $(f_{\xi})\in Z^1(\prod_{j=1}^nH_j,{\rm GL}_NK(x_1,\dots,
x_n))$. Set $\sigma=(\sigma_1,1,\dots,1)\in H_1$ and $\tau=(1,\tau_2,
\dots,\tau_n)\in\prod_{j=2}^nH_j$. By induction assumption 
there exist some $f_1(x),f_2(x)\in{\rm GL}_NK(x)$ such that
$h_{\sigma}(x):=f_1(x)f_{\sigma}(x)f_1(\sigma x)^{-1}\in
{\rm GL}_NK(x_2,\dots,x_n)$ and 
$h_{\tau}(x):=f_2(x)f_{\tau}(x)f_2(\tau x)^{-1}\in{\rm GL}_NK(x_1)$. 

As the sub-semigroups $H_1$ and $\prod_{j=2}^nH_j$ commute, 
$f_{\sigma}(x)f_{\tau}(\sigma x)=f_{\tau}(x)f_{\sigma}(\tau x)$, 
which is equivalent to $$f_1(x)^{-1}h_{\sigma}(x)f_1(\sigma x)
f_2(\sigma x)^{-1}h_{\tau}(\sigma x)f_2(\sigma\tau x)=f_2(x)^{-1}
h_{\tau}(x)f_2(\tau x)f_1(\tau x)^{-1}h_{\sigma}(\tau x)f_1
(\sigma\tau x).$$ 

Set $h(x):=f_1(x)f_2(x)^{-1}$. Then $$[h(x)^{-1}h_{\sigma}(x)h(\sigma x)]
h_{\tau}(\sigma x)=h_{\tau}(x)[h(\tau x)^{-1}h_{\sigma}(\tau x)
h(\sigma\tau x)].$$ Looking at $\tau$'s of type $x_j\mapsto 
x_j^{\ell_j^q}$ and powers of $x_j$ for $2\le j\le n$, we see that 
$$h(x)^{-1}h_{\sigma}(x)h(\sigma x)\in\bigcap_{q\ge 1}{\rm GL}_N
K(x_1,x_2^{\ell_2^q},\dots,x_n^{\ell_n^q})={\rm GL}_NK(x_1),$$ and 
therefore, $$h_{\sigma}(x)[h(\sigma x)h_{\tau}(\sigma x)
h(\sigma x)^{-1}]=[h(x)h_{\tau}(x)h(x)^{-1}]h_{\sigma}(x).$$ 
Looking at $\sigma$'s of type $x_1\mapsto x_1^{\ell_1^q}$ and 
powers of $x_1$, we see that $h(x)h_{\tau}(x)h(x)^{-1}\in
\bigcap_{q\ge 1}{\rm GL}_NK(x_1^{\ell_1^q},x_2,\dots,x_n)=
{\rm GL}_NK(x_2,\dots,x_n)$, and therefore, \begin{equation} \label{odn-iz}
[h_{\sigma}(x),h(x)h_{\tau}(x)h(x)^{-1}]=1. \end{equation}

Fix some $b_2,\dots,b_n\in k$ such that $p(x):=h(x_1,b_2,\dots,b_n)
\in{\rm GL}_NK(x_1)$ is well-defined. Then $h(x)^{-1}h_{\sigma}(x)
h(\sigma x)=p(x)^{-1}h_{\sigma}(x_1,b_2,\dots,b_n)p(\sigma x)$, 
or equivalently, \begin{multline*} q(x):=h(x)p(x)^{-1}=
h_{\sigma}(x)h(\sigma x)p(\sigma x)^{-1}
h_{\sigma}(x_1,b_2,\dots,b_n)^{-1} \\ 
\in\bigcap_{q\ge 1}{\rm GL}_NK(x_1^{\ell_1^q},x_2,\dots,x_n)
={\rm GL}_NK(x_2,\dots,x_n),\end{multline*} 
and in particular, (\ref{odn-iz}) becomes 
$[q(x)^{-1}h_{\sigma}(x)q(x),p(x)h_{\tau}(x)p(x)^{-1}]=1$. 

Set $f(x):=q(x)^{-1}f_1(x)=p(x)f_2(x)$. Then \begin{gather*}r_{\sigma}(x)
:=f(x)f_{\sigma}(x)f(\sigma x)^{-1}=q(x)^{-1}h_{\sigma}(x)q(\sigma x)
=q(x)^{-1}h_{\sigma}(x)q(x)\in{\rm GL}_NK(x_2,\dots,x_n),\\
r_{\tau}(x):=f(x)f_{\tau}(x)f(\tau x)^{-1}=p(x)h_{\tau}(x)p(\tau x)^{-1}
=p(x)h_{\tau}(x)p(x)^{-1}\in{\rm GL}_NK(x_1) \end{gather*} 
and $[r_{\sigma}(x),r_{\tau}(x)]=1$. 

Then the commutation condition $r_{\sigma}(x)r_{\tau}(\sigma x)
=r_{\tau}(x)r_{\sigma}(\tau x)$ implies that 
$$r_{\tau}(x)^{-1}r_{\tau}(\sigma x)=r_{\sigma}(x)^{-1}
r_{\sigma}(\tau x)\in{\rm GL}_NK(x_1)\cap{\rm GL}_NK(x_2,\dots,x_n)
={\rm GL}_NK,$$ so $r_{\sigma}(x)\in\bigcap_{q\ge 1}{\rm GL}_N
K(x_1,x_2^{\ell_2^q},\dots,x_n^{\ell_n^q})={\rm GL}_NK(x_1)$ and 
$$r_{\tau}(x)\in\bigcap_{q\ge 1}{\rm GL}_NK(x_1^{\ell_1^q},x_2,
\dots,x_n)={\rm GL}_NK(x_2,\dots,x_n),$$ and thus, 
$(r_{\xi}(x))\in{\rm Hom}(\prod_{j=1}^nH_j,{\rm GL}_NK)$ 
and it projects onto the class of $(f_{\xi}(x))$. \qed

\begin{lemma} \label{3-3} Let $(f_{\sigma})$ be a 1-cocycle 
on $H:=k^{\times}\rtimes{\mathbb Z}_{\neq 0}$ with values in 
${\rm GL}_Nk(x)$. Suppose that $f_{\xi}(x)$ is regular at 0, 
where $\xi x=x^{\ell}$ for some $\ell\ge 2$, and 
$\mu_{\ell^{\infty}}\subset k$. Then the $L$-semi-linear 
$H$-representation corresponding to $(f_{\sigma})$ is induced 
by a $k$-linear representation of $H$. \end{lemma}
{\it Proof.} In a standard manner, we embed $k(x)$ into $k((x))$. 
Let $f_{\xi}^{-1}=A+B$, where $A\in{\rm GL}_Nk$ 
and $B\in x\mathfrak{gl}_Nk[[x]]$. 

As $A^{-s-1}f_{\xi^{s+1}}^{-1}=A^{-s-1}(A+\xi^sB)f_{\xi^s}^{-1}
=A^{-s}f_{\xi^s}^{-1}+A^{-s-1}\cdot\xi^sB\cdot f_{\xi^s}^{-1}$, 
there is the limit $f:=\lim\limits_{s\to\infty}A^{-s}f_{\xi^s}^{-1}
\in{\bf 1}_N+x\mathfrak{gl}_Nk[[x]]$. 

As explained in \S\ref{contra}, the group ${\rm GL}_NL$ maps bijectively 
onto $Z^1(\langle\xi,\mu_{\ell^{\infty}}\rangle,{\rm GL}_Nk(x))$. 
If $\xi'x=\lambda x^q$ for some $\lambda\in k^{\times}$ 
and $q\in{\mathbb Z}_{\neq 0}$ then $\zeta\xi'=\xi'\zeta^q$ 
for any $\zeta\in\mu_{\ell^{\infty}}$, so 
$f_{\zeta}\cdot\zeta f_{\xi'}=f_{\xi'}\cdot\xi'f_{\zeta^q}$. 
Using expressions for $f_{\zeta}$ we get $f^{-1}\cdot\zeta
(f\cdot f_{\xi'})=f_{\xi'}\cdot\xi'f^{-1}\cdot\xi'\zeta^q f$, 
or equivalently, $f\cdot f_{\xi'}\cdot\xi'f^{-1}=\zeta
(f\cdot f_{\xi'}\cdot\xi'f^{-1})$, and therefore, $f\cdot 
f_{\xi'}\cdot\xi'f^{-1}\in{\rm GL}_Nk((x))^{\mu^n_{\ell^{\infty}}}
={\rm GL}_Nk$, i.e., $f_{\xi'}=f^{-1}\cdot g_{\xi'}\cdot\xi'f$ 
for any $\xi'$ and some 1-cocycle $g_{\xi'}\in{\rm GL}_Nk$. 
%(In particular, $g_{\xi'}={\bf 1}_N$ if $f_{\xi'}(0)={\bf 1}_N$.) 

We may suppose that $k={\mathbb C}$. Then $A^{-s-1}\cdot\xi^sB\cdot 
f_{\xi^s}^{-1}$ can be bounded by 
$c\cdot(N^2\cdot\|A\|\cdot\|A^{-1}\|)^s\cdot|x|^{\ell^s}$, so $f$ is 
holomorphic in a neighbourhood of 0. As 
$f=g_{\xi'}\cdot\xi'f\cdot f_{\xi'}^{-1}$, taking smaller and 
smaller $\lambda\in{\mathbb Q}$, we see that $f$ extends to 
a meromorphic function on ${\mathbb C}$. The involution 
$x\mapsto x^{-1}$ shows that $f$ is rational, and thus, 
$(g_{\sigma})$ belongs to the class of $(f_{\sigma})$. \qed 

\begin{lemma} \label{unip-rep} Let ${\mathbb N}$ act on 
$(k^{\times})^n$ by $\ell(\lambda)=\lambda^{\ell}$ for any 
$\lambda\in(k^{\times})^n$. Then the restriction to $(k^{\times})^n$ 
of any non-degenerate $k$-linear representation $\rho$ of 
$H:=(k^{\times})^n\rtimes{\mathbb N}$ of finite degree $N$ is 
unipotent. (In particular, trivial on $\mu_{\infty}^n$.) \end{lemma}
{\it Proof.} One has $\rho(\lambda)=\rho(\ell)\cdot\rho(\lambda^{\ell})
\cdot\rho(\ell)^{-1}$, so the raising to the $\ell$-th power induces 
a permutation of the eigenvalues of $\rho(\lambda)$. In particular, 
the raising to the $\ell^{N!}$-th power is the identity map of the 
set of eigenvalues of $\rho(\lambda)$, i.e., the eigenvalues of 
$\rho(\lambda)$ form a subset of $\mu_{\ell^{N!}-1}$ for any 
$\ell\ge 2$. Then the eigenvalues of $\rho(\lambda)$ belong to 
the set $\mu_{\ell^{N!}-1}\cap\mu_{(\ell^{N!}-1)^{N!}-1}=\{1\}$, 
and thus, $\rho(\lambda)$ is unipotent. As $\rho(\lambda)$ 
is diagonizable for any $\lambda\in\mu_{\infty}^n$, 
we get $\rho(\lambda)=1$ if $\lambda\in\mu_{\infty}^n$. \qed 

\begin{proposition} \label{unip-tor} $H^1(({\mathbb Z}_{\neq 0}
\ltimes k^{\times})^n,{\rm GL}_Nk)
\stackrel{\sim}{\longrightarrow}H^1(({\mathbb Z}_{\neq 0}\ltimes 
k^{\times})^n,{\rm GL}_Nk(x_1,\dots,x_n))$, if $k$ contains all 
$\ell$-primary roots of unity for some $\ell\ge 2$. \end{proposition}
{\it Proof.} By Theorem \ref{main-loc}, for any extension $K$ 
of $k$ any class in $H^1({\mathbb Z}_{\neq 0}\ltimes k^{\times},
{\rm GL}_NK(t))$ admits a representative with $f_{\xi}$ regular 
at 0, where $\xi t=t^{\ell}$ for some $\ell\ge 2$. 

Then, by Lemma \ref{3-3}, $H^1({\mathbb Z}_{\neq 0}\ltimes 
k^{\times},{\rm GL}_NK)\stackrel{\sim}{\longrightarrow}
H^1({\mathbb Z}_{\neq 0}\ltimes k^{\times},{\rm GL}_NK(t))$. 
By Lemma \ref{red-to-1}, $H^1(({\mathbb Z}_{\neq 0}\ltimes 
k^{\times})^n,{\rm GL}_Nk)$ surjects onto 
$H^1(({\mathbb Z}_{\neq 0}\ltimes k^{\times})^n,
{\rm GL}_Nk(x_1,\dots,x_n))$. \qed 

\begin{corollary} Restriction to 
${\rm Mat}_{n\times n}^{\det\neq 0}{\mathbb Z}\ltimes T$ and 
inclusion ${\rm GL}_Nk\subset{\rm GL}_NL$ induce a natural bijection 
$H^1({\rm GL}_n{\mathbb Q}\ltimes(T\otimes{\mathbb Q}),{\rm GL}_Nk)
\stackrel{\sim}{\longrightarrow}
H^1({\rm Mat}_{n\times n}^{\det\neq 0}{\mathbb Z}\ltimes T,
{\rm GL}_NL)$, if $k$ contains all $\ell$-primary roots 
of unity for some $\ell\ge 2$, where $L:=k(x_1,\dots,x_n)$ and 
$T=(k^{\times})^n$. \end{corollary}
{\it Proof.}  By Proposition \ref{unip-tor}, any class in 
$H^1({\rm Mat}_{n\times n}^{\det\neq 0}{\mathbb Z}\ltimes T,
{\rm GL}_NL)$ is represented by a cocycle with constant 
restriction to ${\mathbb Z}_{\neq 0}^n\ltimes T$. For any 
$\lambda\in{\rm Mat}_{n\times n}^{\det\neq 0}{\mathbb Z}$ and any 
$\mu\in T_{{\rm tors}}$ one has 
$\lambda\circ\mu^{\lambda}=\mu\circ\lambda$, so by Lemma \ref{unip-rep}, 
$f_{\lambda}(x)=f_{\lambda}(x)f_{\mu^{\lambda}}(x^{\lambda})=
f_{\mu}(x)f_{\lambda}(\mu\cdot x)=f_{\lambda}(\mu\cdot x)$, i.e., 
$f_{\lambda}(x)\in{\rm GL}_Nk(x_1,\dots,x_n)^{T_{{\rm tors}}}
={\rm GL}_Nk$. \qed 

\section{Main theorem} \label{techno}
\begin{lemma} \label{prel} Let $h(t)\in{\rm GL}_Nk(t)$ be 
a function in one variable such that $h(x^{-1})=h(x)^{-1}$ 
and $h(x)h(y)=h(xy-x+1)h(xy(xy-x+1)^{-1})$. Then 
$h$ is regular on ${\mathbb G}_m-\{1\}$; 

the values of $h$ generate a closed algebraic connected subgroup 
$H$ of ${\rm GL}_Nk$; 

the unipotent radical $H_u$ of $H$ is commutative, coincides with 
the commutant of $H$ and with $\{h(t_1)\cdots h(t_m)~|~m\ge 2,
~t_1,\dots,t_m\in k^{\times}-\{1\},~t_1\cdots t_m=1\}$. \end{lemma} 
{\it Proof.} Let $\Sigma$ be the set of poles of $h(x)^{\pm 1}$, 
$a\in k^{\times}-\{1\}$ and $b\in k^{\times}$ be outside of the finite 
set $\Sigma\cup(\frac{1}{a}\Sigma-\frac{1}{a}+1)\cup\frac{a-1}{a}
(1-\Sigma^{-1})^{-1}$. Then $h(x)=h((b-1)x+1)
h(bx((b-1)x+1)^{-1})h(b)^{-1}$ is regular at $a$, so 
$\Sigma\subseteq\{0,1,\infty\}$. 

Let $x_1,\dots,x_m$ be variables. Set 
$X_j:=x_1\cdots x_{j-1}(x_j-1)$ for all $1\le j\le m$. 
Let $\sigma_i$ be the automorphism of $L:=k(x_1,\dots,x_m)
=k(X_1,\dots,X_m)$ over $k$ such that $$\sigma_ix_j:=\left\{
\begin{array}{ll}x_j & \mbox{if $i\le j-2$, or $i\ge j+1$} \\ 
x_j(x_{j+1}-1)+1 & \mbox{if $i=j$} \\ 
x_jx_{j+1}(x_jx_{j+1}-x_j+1)^{-1} & \mbox{if $i=j-1$} 
\end{array}\right.$$
Then $\sigma_i(h(x_1)\cdots h(x_m))=h(x_1)\cdots h(x_m)$ for all 
$1\le i\le m-1$. Denote by ${\mathfrak S}$ the group of automorphisms 
of $L$ generated by $\sigma_1,\dots,\sigma_{m-1}$. As $$\sigma_iX_j:=
\left\{\begin{array}{ll} x_1\cdots x_{j-2}(x_{j-1}x_j-x_{j-1}+1)
\left(\frac{x_{j-1}x_j}{x_{j-1}x_j-x_{j-1}+1}-1\right)=X_{j-1} & 
\mbox{if $i=j$} \\ 
x_1\cdots x_{j-1}x_j(x_{j+1}-1)=X_{j+1} & \mbox{if $i=j-1$} \\ 
X_j & \mbox{otherwise} 
\end{array}\right.$$ 
${\mathfrak S}$ is isomorphic to the symmetric group ${\mathfrak S}_m$, 
and $h(x_1)\cdots h(x_m)\in{\rm GL}_NL^{{\mathfrak S}}$. 

Let $\sigma\in{\mathfrak S}$ be such that $\sigma X_j=X_{m-j+1}$ 
for all $1\le j\le m$. Then $\sigma(X_j+X_{j-1}+\cdots+X_1)=
X_{m-j+1}+X_{m-j+2}+\cdots+X_m=x_1\cdots x_m-x_1\cdots x_{m-j}$. 
As $x_j=\frac{X_j+X_{j-1}+\cdots+X_1+1}{X_{j-1}+\cdots+X_1+1}$, 
we get $\sigma x_j=\frac{x_1\cdots x_m-x_1\cdots x_{m-j}+1}
{x_1\cdots x_m-x_1\cdots x_{m-j+1}+1}$, and thus, 
\begin{multline*} h(x_1)\cdots h(x_m)=h(x_1\cdots x_m-x_1\cdots 
x_{m-1}+1)h\left(\frac{x_1\cdots x_m-x_1\cdots x_{m-2}+1}
{x_1\cdots x_m-x_1\cdots x_{m-1}+1}\right) \\ 
\cdots h\left(\frac{x_1\cdots x_m-x_1+1}{x_1\cdots x_m-x_1x_2+1}\right)
h\left(\frac{x_1\cdots x_m}{x_1\cdots x_m-x_1+1}\right).\end{multline*} 

Then for any $t_1,\dots,t_m\in k^{\times}-\{1\}$ such that 
$t_1\cdots t_m=-1$ we get $$h(t_1)\cdots h(t_m)=
h(t_m^{-1})\cdots h(t_1^{-1})=(h(t_1)\cdots h(t_m))^{-1},$$ 
or equivalently, $(h(t_1)\cdots h(t_m))^2=1$. 

Let $H$ be the subgroup in ${\rm GL}_Nk$ generated by $h(t)$
for all $t\in k^{\times}-\{1\}$, and $S=\{h(t_1)\cdots h(t_m)
~|~m\ge 2;t_1,\dots,t_m\in k^{\times}-\{1\};t_1\cdots t_m=1\}$. 
Clearly, $S$ is a subgroup of $H$. One has $h(-1)Ah(-1)=A^{-1}$ 
for any $A\in S$, so $h(-1)ABh(-1)=A^{-1}B^{-1}=B^{-1}A^{-1}$ 
for any $A,B\in S$, which means that $S$ is abelian. 

Let $({\mathbb G}_m-\{1\})^M\stackrel{\pi_M}{\longrightarrow}H$ 
sends $(t_1,\dots,t_M)$ to $h(t_1)\cdots h(t_M)$ and 
$P_M:=\{(t_1,\dots,t_M)\in({\mathbb G}_m-\{1\})^M~|~t_1\cdots t_M=1\}$. 

Then $\pi_{M+2}(t_1,\dots,t_M,t,t^{-1})=\pi_M(t_1,\dots,t_M)$ 
and $\pi_{M+3}(t_1,\dots,t_M,t,1+\zeta-\zeta t^{-1},
((1+\zeta)t-\zeta)^{-1})=\pi_M(t_1,\dots,t_M)$, 
so ${\rm Im}(\pi_{M-3})\subseteq{\rm Im}(\pi_M)\supseteq
{\rm Im}(\pi_{M-2})$ and ${\rm Im}(\pi_{M-3}|_{P_{M-3}})\subseteq
{\rm Im}(\pi_M|_{P_M})\supseteq{\rm Im}(\pi_{M-2}|_{P_{M-2}})$. As 
$({\mathbb G}_m-\{1\})^M$ and its subset $P_M$ are irreducible, this 
implies that $H$ and $S$ are closed algebraic connected subgroups 
of ${\rm GL}_Nk$. Clearly, $[H,H]\subseteq S$, so $H$ is 
solvable. 

As $S$ is connected, one has $S=S^2$ (the set of the squares of 
the elements of $S$), so, in view of $[H,H]\supseteq[h(-1),S]=S^2$, 
we get $S=[H,H]$. By Lie--Kolchin theorem, we may conjugate $H$
by an element of ${\rm GL}_Nk$ so that it becomes upper 
triangular, and $S$ is contained in the unipotent radical of $H$. 

If $N=1$ then $S=[H,H]=\{1\}$, so $h(x)h(y)=h(xy)$, and thus, $h$ 
is a cocharacter. This implies that for any $N\ge 1$ the diagonal 
of $h$ is a cocharacter. 

As the unipotent radical $H_u$ of $H$ coincides 
with the set of unipotent elements of $H$, one has 
$$H_u=\{h(t_1)\cdots h(t_m)~|~(t_1\cdots t_m)^{(m_1,\dots,m_N)}=1\}
=S\cup\bigcup_{\zeta\in\mu_{(m_1,\dots,m_N)}-\{1\}}h(\zeta)S,$$ 
where the diagonal of $h(t)$ is ${\rm diag}(t^{m_1},\dots,t^{m_N})$. 
As $H_u$ is connected, one has $H_u=S$. \qed 

\begin{lemma} \label{10half} Keeping notations of Lemma \ref{prel}, 
let $f\in{\rm GL}_Nk$ be a unipotent matrix normalizing $H$. 
Suppose that $h(x)h(2-x^{-1})=f\cdot h(2x-1)\cdot f^{-1}$. 
Then $h(x)h(y)=h(xy)$. \end{lemma} 
{\it Proof.} Zariski closure of the subgroup generated by $f$ is 
connected, as well as any other closed abelian unipotent subgroup. 
Let $\tilde{H}$ be the Zariski closed subgroup generated by $f$ 
and $H$. Then $\tilde{H}$ is solvable and connected, so we may 
assume that it is upper triangular. 

By induction on $N$, the case $N=1$ being trivial, 
we check that there is a diagonal matrix 
in the ${\rm GL}_Nk$-conjugacy class of $h(t)$. 

Let $N\ge 2$. One has $h(t)=\left(\begin{array}{cc} 
\Lambda(t)^{-1} & g(t) \\ 0 & M(t)\end{array}\right)\in
\left(\begin{array}{cc} {\rm GL}_{N-1}k(t) & k(t)^{N-1} \\ 
0 & k(t)^{\times}\end{array}\right)$, where 
an upper triangular $\Lambda(t)^{-1}\in{\rm GL}_{N-1}k(t)$ 
and $M(t)\in k(t)^{\times}$ satisfy the same conditions as 
$h(t)$, so we may suppose that $\Lambda(t)$ is diagonal. 

The condition $h(t)^{-1}=h(t^{-1})$ gives $g(t^{-1})=-\Lambda(t)
g(t)M(t)^{-1}$. Set $v(t):=\Lambda(t)g(t)$. Then 
$v(t)=-g(t^{-1})M(t)=-\Lambda(t)v(t^{-1})M(t)$. 

One has $$h(x)h(y)=\left(\begin{array}{cc} 
\Lambda(xy)^{-1} & \Lambda(x)^{-1}g(y)+g(x)M(y) \\ 
0 & M(xy)\end{array}\right),$$ so $\Lambda(x)^{-1}g(y)+g(x)M(y)=
\Lambda(xy-x+1)^{-1}g(xy(xy-x+1)^{-1})+g(xy-x+1)M(xy(xy-x+1)^{-1})$. 
Multiplying both sides by $\Lambda(x)$ on the left and by $M(y)^{-1}$ 
on the right, we get $v(x)-v(y^{-1})=
\Lambda(y)^{-1}v(xy(xy-x+1)^{-1})M(y)^{-1}-\Lambda(x)v((xy-x+1)^{-1})M(x)$, 
or equivalently, \begin{equation}\label{razn} 
v(x)-v(y)=\Lambda(y)v(X)M(y)-\Lambda(x)v(Y)M(x),\end{equation} 
where we set $X=x(x-xy+y)^{-1}$ and $Y=y(x-xy+y)^{-1}$. 

Taking the partial derivative $\partial/\partial x$ of both sides 
of (\ref{razn}) gives \begin{multline} \label{per-der} v'(x)=
\left(\frac{1}{x-xy+y}-\frac{x(1-y)}{(x-xy+y)^2}\right)
\Lambda(y)v'(X)M(y)-\Lambda'(x)v(Y)M(x) \\ 
-\Lambda(x)v(Y)M'(x)+\frac{y(1-y)}{(x-xy+y)^2}
\Lambda(x)v'(Y)M(x).\end{multline} 

Taking further the partial derivative $\partial/\partial y$ of 
both sides of (\ref{per-der}) gives \begin{multline*} 
\frac{y}{(x-xy+y)^2}\Lambda'(y)v'(X)M(y)+\left(\frac{1}{(x-xy+y)^2}
-\frac{2(1-x)y}{(x-xy+y)^3}\right)\Lambda(y)v'(X)M(y)\\ 
+\frac{x(x-1)y}{(x-xy+y)^4}\Lambda(y)v''(X)M(y)+
\frac{y}{(x-xy+y)^2}\Lambda(y)v'(X)M'(y) \\ 
+\left(\frac{y(1-x)}{(x-xy+y)^2}-\frac{1}{x-xy+y}\right)
\Lambda'(x)v'(Y)M(x)-\frac{x}{(x-xy+y)^2}
\Lambda(x)v'(Y)M'(x) \\ 
+\left(\frac{1-2y}{(x-xy+y)^2}-\frac{2(1-x)y(1-y)}{(x-xy+y)^3}
\right)\Lambda(x)v'(Y)M(x) \\ 
+\frac{xy(1-y)}{(x-xy+y)^4}\Lambda(x)v''(Y)M(x)=0. \end{multline*} 

As $X(X-1)=\frac{x(x-1)y}{(x-xy+y)^2}$, 
$Y(Y-1)=\frac{xy(y-1)}{(x-xy+y)^2}$, $2X-1=\frac{x-y+xy}{x-xy+y}$ 
and $1-2Y=\frac{x-y-xy}{x-xy+y}$, multiplying by $(x-xy+y)^2$, 
we get \begin{multline*} X(X-1)\Lambda(y)v''(X)M(y)
-Y(Y-1)\Lambda(x)v''(Y)M(x)\\ +y\Lambda'(y)v'(X)M(y)
-x\Lambda'(x)v'(Y)M(x)+(2X-1)\Lambda(y)v'(X)M(y)\\
-(2Y-1)\Lambda(x)v'(Y)M(x)+y\Lambda(y)v'(X)M'(y)
-x\Lambda(x)v'(Y)M'(x)=0. \end{multline*} 

As $x/y=X/Y$, multiplying by $\Lambda(X)^{-1}\Lambda(y)^{-1}$ 
on the left and by $M(y)^{-1}M(X)^{-1}$ on the right, and using 
$y\Lambda'(y)\Lambda(y)^{-1}=\Lambda'(1)$ and $yM'(y)M(y)^{-1}=M'(1)$, 
we see that the function \begin{multline} \label{con-fun}
X(X-1)\Lambda(X)^{-1}v''(X)M(X)^{-1}
+\Lambda(X)^{-1}\Lambda'(1)v'(X)M(X)^{-1}\\ 
+(2X-1)\Lambda(X)^{-1}v'(X)M(X)^{-1}
+\Lambda(X)^{-1}v'(X)M'(1)M(X)^{-1} \end{multline} 
is invariant under $X\leftrightarrow Y$, i.e., it is constant. 

Let $v'(X)=X^{-1}(X-1)^{-1}\Lambda(X(X-1)^{-1})\varphi(X)
M(X(X-1)^{-1})$. Then \begin{multline*} v''(X)=-\left(\frac{1}{X}+
\frac{1}{X-1}\right)v'(X)-\frac{1}{(X-1)^2}\Lambda'\left(\frac{X}{X-1}
\right)\Lambda\left(\frac{X-1}{X}\right)v'(X)\\ -\frac{1}{(X-1)^2}v'(X)
M\left(\frac{X-1}{X}\right)M'\left(\frac{X}{X-1}\right)
+\frac{1}{X(X-1)}\Lambda\left(\frac{X}{X-1}\right)\varphi'(X)
M\left(\frac{X}{X-1}\right),\end{multline*} 
and therefore, the function (\ref{con-fun}) coincides with 
$$\Lambda(X-1)^{-1}\varphi'(X)M(X-1)^{-1}=A\in k^{N-1}.$$ 
This implies that $\varphi'(X)=\Lambda(X-1)AM(X-1)$ is a 
Laurent polynomial in $X-1$. It follows from the rationallity of 
$\varphi$ that $\varphi(X)=(X-1)\Lambda(X-1)BM(X-1)+C$ for some 
$B,C\in k^{N-1}$, and therefore, $v'(X)=X^{-1}\Lambda(X)BM(X)+
X^{-1}(X-1)^{-1}\Lambda(X(X-1)^{-1})CM(X(X-1)^{-1})$. 

It follows from the rationallity of $v$ that $$v(X)=\Lambda(X)DM(X)+
\Lambda(X(X-1)^{-1})EM(X(X-1)^{-1})+F$$ for some $D,E,F\in k^{N-1}$. 

The condition $v(t^{-1})=-\Lambda(t)^{-1}v(t)M(t)^{-1}$ is equivalent 
to $D+\Lambda(X-1)^{-1}EM(X-1)^{-1}+\Lambda(X)^{-1}FM(X)^{-1}+
\Lambda(X)^{-1}DM(X)^{-1}+\Lambda(1-X)^{-1}EM(1-X)^{-1}+F=0$. 
Rewriting it as $(D+F)+\Lambda(X)(D+F)M(X)+\Lambda(X(X-1)^{-1})
(E+\Lambda(-1)EM(-1))M(X(X-1)^{-1})=0$, we see that adding to 
$(D,E,F)$ of a multiple of $(D+F,E+\Lambda(-1)EM(-1),D+F)$ 
does not change $v$, so we may assume that $F=-D$ and 
$E+\Lambda(-1)EM(-1)=0$, and thus, $g(t)=DM(t)+\Lambda(t-1)^{-1}
EM(t(t-1)^{-1})-\Lambda(t)^{-1}D$. 

As \begin{multline*} \left(\begin{array}{cc} {\bf 1} & -D \\ 
0 & 1 \end{array}\right)
\left(\begin{array}{cc} \Lambda^{-1} & g \\ 0 & M \end{array}\right) 
\left(\begin{array}{cc} {\bf 1} & D \\ 0 & 1 \end{array}\right)
=\left(\begin{array}{cc} \Lambda^{-1} & g-DM+\Lambda^{-1}D \\ 
0 & M \end{array}\right)\\ =\left(\begin{array}{cc} 
\Lambda(t)^{-1} & \Lambda(t-1)^{-1}EM\left(\frac{t}{t-1}\right) \\ 
0 & M \end{array}\right),\end{multline*} we may suppose 
that $D=0$ and $E_j=0$ if $m_j-m_N=1$ (since then 
$(\Lambda(t-1)^{-1}EM(\frac{t}{t-1}))_j=(\Lambda(t)^{-1}E-EM(t))_j$). 

It is easy to see that $h(x)h(y)=h(xy-x+1)h(xy(xy-x+1)^{-1})$. 
The condition $h(t^{-1})=h(t)^{-1}$ is equivalent to 
$\Lambda(t(1-t)^{-1})EM(1-t)^{-1}=
-\Lambda(t(t-1)^{-1})EM(t-1)^{-1}$, 
which is the same as $E=-\Lambda(-1)EM(-1)$, and thus, 
$E_j=0$ if $m_j-m_N$ is even. 

Now suppose that for some unipotent $f=\left(\begin{array}{cc} 
A & B \\ 0 & 1 \end{array}\right)\in{\rm GL}_Nk$ one has 
$$f\left(\begin{array}{cc} \Lambda(x)^{-1} & 
\Lambda(x-1)^{-1}EM\left(\frac{x}{x-1}\right) \\ 0 & M(x) 
\end{array}\right) =\left(\begin{array}{cc} \Lambda(x)^{-1} & 
2\Lambda\left(\frac{2}{x-1}\right)EM\left(\frac{2x}{x-1}\right) \\ 
0 & M(x) \end{array}\right)f.$$
Then $A\Lambda(x)^{-1}=\Lambda(x)^{-1}A$, which is equivalent 
to $A_{ij}=0$ if $m_i\neq m_j$, and $$A\Lambda(x-1)^{-1}
EM(\frac{x}{x-1})+BM(x)=\Lambda(x)^{-1}B+
2\Lambda(\frac{2}{x-1})EM(\frac{2x}{x-1})D.$$

These conditions imply that $B_i(1-x^{m_i-m_N})=(2^{m_N-m_i+1}E_i
-\sum\limits_{j:m_i=m_j}A_{ij}E_j)(x-1)^{m_i-m_N}$, and thus, 
$B_i=0$ and $2^{m_N-m_i+1}E_i=\sum\limits_{j:m_i=m_j}
A_{ij}E_j$ if $m_i-m_N\not\in\{0,1\}$. 

This implies that $2^{m_N-m_i+1}E_i=(AE)_i$ for all $i$, 
and therefore, $\Lambda(2)^{-1}AE=2^{m_N+1}E$, or equivalently, 
$(A-1)E=(2^{m_N+1}\Lambda(2)-1)E$. Then $0=(A-1)^{N-1}E=
(2^{m_N+1}\Lambda(2)-1)^{N-1}E$, and thus, 
$(\Lambda(2)^{-1}-2^{m_N+1})E=0$, and therefore, 
$(m_i-m_N-1)E_i=0$ for all $i$. This means that $E=0$, i.e., 
$h(t)$ is diagonal. \qed 

\begin{lemma} \label{restr-uni-rad} Let $\mu_{\ell^{\infty}}\subset k$ 
for some $\ell\ge 2$ and $L:=k({\mathbb P}^n_k)$. \begin{enumerate} 
\item Let $P$ be the stabilizer of a hyperplane $H$ in 
${\mathbb P}^n_k$, and $U$ the unipotent radical of $P$. 
Let $V$ be an $L$-semi-linear $P$-representation. Suppose that 
the restriction of $V$ to a maximal torus $T\subset P$ is induced 
by a unipotent representation of degree $N$. Then there is a 
$U$-invariant\footnote{more generally, invariant under the normalizer 
in $P$ of a $k_{\infty}$-lattice in $U$.} $k$-lattice $V_0$ in $V$ with 
unipotent action of $U$ such that 
${\mathcal O}_{{\mathbb P}^n_k}({\mathbb P}^n_k-H)\otimes_kV_0$ 
is canonical and $P$-invariant. 
\item \label{equi-le} If, moreover, $V$ is an $L$-semi-linear 
$G$-representation, where $G={\rm Aut}({\mathbb P}^n_k/k)$ with 
$n\ge 2$, then $V$ is the generic fibre of a coherent $G$-sheaf 
on ${\mathbb P}^n_k$. \end{enumerate}
%representation $P/U\stackrel{\rho}{\longrightarrow}{\rm GL}_Nk$. 
\end{lemma}
{\it Proof.} \begin{enumerate} \item We shall consider $U$ as 
a $k$-vector space and identify the action of ${\rm GL}_k(U)$ 
with the adjoint action of $P/U\cong{\rm GL}_nk$. 

Let $(f_{\sigma})\in H^1(P,{\rm GL}_NL)$ be the class of $V$. 
We may suppose that $f_{\lambda}\in{\rm GL}_Nk$ is 
unipotent for any $\lambda\in T$, and in particular, $f_{\lambda}=1$ 
for any $\lambda\in T_{{\rm tors}}$. Then $f_{\lambda u}(x)=
f_{\lambda}^{-1}f_u(\lambda^{-1}x)f_{\lambda}$ for any $u\in U$ and 
$\lambda\in T$. Choose a splitting $T=({\mathbb G}_m)^n$ (which is 
equivalent to ordering of $(n+1)$ points of ${\mathbb P}^n_k$ fixed 
by $T$) and for each $1\le j\le n$ choose some non-zero $u_j\in U$ 
fixed by $T_j:=\mu_{\ell^{\infty}}^{j-1}\times\{1\}\times
\mu_{\ell^{\infty}}^{n-j}$. 
Then $f_{u_j}(x)\in{\rm GL}_NL^{T_j}={\rm GL}_Nk(x_j)$, where 
$\tau x_j=\tau_jx_j$ for any $\tau=(\tau_1,\dots,\tau_n)\in T$, 
which implies that $f_{u_i}(x)$ and $f_{u_j}(y)$ commute 
for any $i\neq j$. 

Set $h_j(t):=f_{u_j}(\frac{u_j}{t-1})\in{\rm GL}_Nk(t)$. Then, as 
$$f_{u_j}(x)f_{\lambda u_j}(x+u_j)=f_{(1+\lambda)u_j}(x)=
f_{\lambda u_j}(x)f_{u_j}(x+\lambda u_j)$$ for any $\lambda\in 
k^{\times}\subseteq T$, we get \begin{multline}\label{usl-koc}
h_j(t)\cdot f_{\lambda}^{-1}h_j(1+\lambda-\lambda\cdot t^{-1})
f_{\lambda}=f_{1+\lambda}^{-1}h_j((1+\lambda)t-\lambda)f_{1+\lambda}\\ 
=f_{\lambda}^{-1}h_j(\lambda(t-1)+1)f_{\lambda}\cdot 
h_j(\frac{(1+\lambda)t-\lambda}{\lambda(t-1)+1})\end{multline} 
if $\lambda\in k^{\times}-\{-1\}$, and $h_j(t^{-1})=h_j(t)^{-1}$ 
if $\lambda=-1$. 

If $\lambda=1$ this gives $h_j(t)h_j(2-t^{-1})=f^{-1}_2h_j(2t-1)f_2$. 
Setting $y:=1+\lambda-\lambda t^{-1}$, we get $h_j(t)h_j(y)=h_j(ty-t+1)
h_j(\frac{ty}{ty-t+1})$ if $\lambda\in\mu_{\ell^{\infty}}$. As 
$\mu_{\ell^{\infty}}$ is Zariski dense in ${\mathbb G}_m$, this 
identity holds for arbitrary $t$ and $y$. Then $h_j(t)$ satisfies 
the conditions of Lemma \ref{10half}, so $h_j(x)h_j(y)=h_j(xy)$. 
As $f_{u_i}(x)$ and $f_{u_j}(y)$ commute, the same holds for 
$h_i(x)$ and $h_j(y)$, so $h_j(t)=
C^{-1}\cdot{\rm diag}(t^{m_{1j}}_j,\dots,t^{m_{Nj}}_j)\cdot C$ 
for some $C\in{\rm GL}_Nk$ and some $m_{ij}\in{\mathbb Z}$. 
This is equivalent to $f_{u_j}(x)=C^{-1}\cdot{\rm diag}
((1+x_j(u_j)/x_j)^{m_{1j}},\dots,(1+x_j(u_j)/x_j)^{m_{Nj}})
\cdot C$. 

One sees from (\ref{usl-koc}) that if $f_{\lambda}$ and $f_{\mu}$ 
commute with $h_j(t)$ for some $\lambda\neq-\mu$ in $k^{\times}$ then 
$f_{\lambda+\mu}$ and $f_{\lambda\mu}$ also commute with $h_j(t)$. 
For any root $\mu$ of $\lambda$ in $k$ the element $f_{\mu}$ 
belongs to the Zariski closure of the subgroup generated by 
$f_{\lambda}$, so $f_{\lambda}$ commutes with $h_j(t)$ for any 
$\lambda$ in the radical closure $k_{\infty}$ of ${\mathbb Q}$ 
in $k$ (defined in \S\ref{int-1}). 

Set $g(x)=\prod_{j=1}^nh_j(x_j)$ and 
$g_{\lambda}(x):=g(x)f_{\lambda}(x)g(\lambda x)^{-1}$ for any 
$\lambda\in P$. Then $g_u(x)=1$ for any $u$ in the $k_{\infty}$-subspace 
$U_0$ of $U$ spanned by $u_1,\dots,u_n$. 

For any $\lambda\in P$ and any $u\in U$ one has $g_{\lambda}(x)
g_{\lambda^{-1}u\lambda}(\lambda x)=g_u(x)g_{\lambda}(x+u)$. 
If $\lambda\in U_0$ we see that $g_u$ is constant. As $U$ is normal 
in $P$, this shows that, for arbitrary $\lambda\in P$, the poles of 
$g_{\lambda}(x)^{\pm 1}$ are $U$-invariant, so $g_{\lambda}(x)
\in{\rm GL}_Nk[x]$. In other words, the submodule 
${\mathcal V}_H:={\mathcal O}_{{\mathbb P}^n_k}
({\mathbb P}^n_k-H)\otimes_kV^{U_0}$ is $P$-invariant. As all 
$k_{\infty}$-lattices $U_0$ form a single $P$-orbit, ${\mathcal V}_H$ 
is independent of the choice of $U_0$. %$g_{\lambda}(x)=g_{\lambda}(x+u)$, 
%In other words, $k[x]\otimes_kV^{\langle u_1,\dots,u_n\rangle}$ 
%is a complete $k[x]$-lattice in $V$. 

For any $\lambda\in P$ and any $u\in U_0$ one has 
$g_{\lambda}(x)^{-1}g_{\lambda}(x+u)=g_{\lambda^{-1}u\lambda}
\in{\rm GL}_Nk$, so one sees that $g_{\lambda}(x)\in{\rm GL}_Nk$ 
for any $\lambda$ in the normalizer of $U_0$ in $P$ and that 
$r_{\lambda}(u):=g_{\lambda}(x)^{-1}g_{\lambda}(x+u)\in{\rm GL}_Nk[u]$. 
Then $g_{\lambda}(x)=g_{\lambda}(0)r_{\lambda}(x)$. Note that 
$$r_{\lambda}(u)r_{\lambda}(v)=g_{\lambda}(x)^{-1}g_{\lambda}(x+u)
g_{\lambda}(x+u)^{-1}g_{\lambda}(x+u+v)=r_{\lambda}(u+v),$$
i.e., $r_{\lambda}:U\longrightarrow{\rm GL}_Nk$ is a polynomial 
homomorphism and $r_{\lambda}=1$ for any $\lambda\in U$. 

In particular, the restriction of $V$ to $U$ is induced by 
a unipotent representation trivial on $U_0$. 
\item As we know from the first part, for any hyperplane $H$ and 
any $k_{\infty}$-lattice $U_0$ in the unipotent radical $U$ of ${\rm Stab}_H$ 
the submodule ${\mathcal V}_H:={\mathcal O}_{{\mathbb P}^n_k}
({\mathbb P}^n_k-H)\otimes_kV^{U_0}$ is canonical and 
${\rm Stab}_H$-invariant. Then the submodule ${\mathcal V}_p:=
{\mathcal O}_{{\mathbb P}^n_k,p}\otimes_kV^{U_0}$ of $V$ is also 
independent of the choice of $U_0$, and the group 
${\rm Stab}_p\cap{\rm Stab}_H$ acts on ${\mathcal V}_p$. 
Let, in notation of the first part, $\tau_u:x\mapsto\frac{x}{1+\langle 
u,x\rangle}$ be an element of the unipotent radical of ${\rm Stab}_p$, 
where $u$ is a $k$-linear functional on $U$. Then 
$g_{\lambda}(x)g_{\lambda^{-1}\tau_u\lambda}(\lambda x)
=g_{\tau_u}(x)g_{\lambda}(\frac{x}{1+\langle u,x\rangle})$ for any 
$\lambda\in{\rm Stab}_p\cap{\rm Stab}_H$. As $g_{\lambda}(x),g_{\lambda}
(\frac{x}{1+\langle u,x\rangle})\in{\rm GL}_N
{\mathcal O}_{{\mathbb P}^n_k,p}$, we see that the singularities 
of $g_{\tau_u}(x)^{\pm 1}$ in the formal neighbourhood of $p$ are 
invariant under the action of the centralizer of $\tau_u$ in 
${\rm Stab}_p\cap{\rm Stab}_H$, i.e., $g_{\tau_u}(x)\in
{\rm GL}_N{\mathcal O}_p[\langle u,x\rangle^{-1}]$. 

We see from  $g_{\tau_u}(x)
g_{\tau_v}(\frac{x}{1+\langle u,x\rangle})=g_{\tau_{u+v}}(x)$ that 
$g_{\tau_u}(x)\in{\rm GL}_N{\mathcal O}_p$ if $n\ge 2$, i.e., that 
${\mathcal V}_p$ is invariant under the action of the unipotent 
radical of ${\rm Stab}_p$, and thus, under the action of ${\rm Stab}_p$ 
itself. As all hyperplanes avoiding a point $p$ form a single 
${\rm Stab}_p$-orbit, the submodule ${\mathcal O}_p\otimes_kV^{U_0}$ 
of $V$ is independent of the choice of the hyperplane, 
so we get a coherent sheaf ${\mathcal V}\subset V$. 
%For any $u\in U_0$ one has $g_u(0)=1$, so rewriting the 1-cocycle 
%condition for $(g_{\lambda})$ as $g_{\lambda}(0)r_{\lambda}(x)
%g_{\mu}(0)r_{\mu}(\lambda x)=g_{\lambda\mu}(0)r_{\lambda\mu}(x)$, 
%we see that $g(0):P\stackrel{\lambda\mapsto g_{\lambda}(0)}
%{-\!\!\!-\!\!\!-\!\!\!\longrightarrow}{\rm GL}_Nk$ is a 
%homomorphism trivial on the subgroup $U_0$. As $P$ acts on $U$ 
%transitively, $g(0)$ is trivial on entire $U$, so $V=L\otimes_kV^U$ 
%and $g_{\lambda}\in{\rm GL}_Nk$ for any $\lambda\in P$. 
\qed \end{enumerate}

\begin{lemma} \label{predp} Let $\mu_{\ell^{\infty}}\subset k$ 
for some $\ell\ge 2$ and $L:=k({\mathbb P}^n_k)$ with $n\ge 2$. 
Let $V$ be an $L$-semi-linear $G$-representation, where $G=
{\rm Aut}({\mathbb P}^n_k/k)\cong{\rm PGL}_{n+1}k$. Suppose 
that the restriction of $V$ to the stabilizer $P$ of a hyperplane 
in ${\mathbb P}^n_k$ is induced by a representation of degree 
$N$ trivial on the unipotent radical $U$ of $P$. Then $V$ 
is isomorphic to the $L$-semi-linear $G$-representation 
corresponding to the generic fibre of a $G$-equivariant 
coherent sheaf on ${\mathbb P}^n_k$ of rank $N$. \end{lemma}
{\it Proof.} There is a representative $(g_{\lambda})\in 
Z^1(G,{\rm GL}_NL)$ of $V$ such that $g_{\lambda}\in{\rm GL}_Nk$ 
for any $\lambda\in P$ and $g_{\lambda}=1$ for any $\lambda\in U$. 
As $G$ is generated by a finite number of conjugates of $U$, the 
${\rm GL}_N$-valued function $g_{\lambda}(x)$ is rational on 
$U\times_kG$, which implies that the coherent sheaf of Lemma 
\ref{restr-uni-rad} (\ref{equi-le}) is $G$-equivariant. \qed 

\begin{theorem} \label{main-th} Let $\mu_{\ell^{\infty}}\subset k$ 
for some $\ell\ge 2$ and $L:=k({\mathbb P}^n_k)$. 
Let $n\ge 2$ and $V$ be a (non-degenerate) finite-dimensional 
$L$-semi-linear ${\rm End}(L/k)$-representation. Then the restriction 
of $V$ to $G:={\rm Aut}({\mathbb P}^n_k/k)\cong{\rm PGL}_{n+1}k$ 
is isomorphic to the $L$-semi-linear $G$-representation 
corresponding to the generic fibre of a coherent $G$-sheaf on 
${\mathbb P}^n_k$, $G$-equivariant if $k=k_{\infty}$. \end{theorem}
{\it Proof.} Fix a maximal torus $T\subset G$ and corresponding 
coordinates $x_1,\dots,x_n$. In fact, $V$ is an $L$-semi-linear 
representation of the semi-group generated by $G$ and 
$\sigma_p:x\mapsto x^p$ for all $p\in{\mathbb Z}_{\neq 0}^n$. 
By Proposition \ref{unip-tor}, the restriction of $V$ to $T\rtimes
{\mathbb Z}_{\neq 0}^n$ is induced by a representation. Then, 
by Lemma \ref{unip-rep}, the restriction of $V$ to $T$ is induced 
by a unipotent representation. Lemma \ref{restr-uni-rad} produces 
a coherent $G$-sheaf on ${\mathbb P}^n_k$, 
and Lemma \ref{predp} implies the rest. \qed 

\vspace{5mm}

{\it Remark.} \label{podfakt} Let $Q$ be a finite-dimensional
$k$-vector space, $G={\rm PGL}(Q)$ and ${\mathcal V}$ be a 
coherent $G$-sheaf on ${\mathbb P}(Q)$. Then there exists a 
$k$-linear representation $W$ of $G$ of finite degree and an 
isomorphism of ${\mathcal V}$ onto a $G$-subquotient of
$W\otimes_k{\mathcal O}_{{\mathbb P}(Q)}$. ({\it Proof.} For $m$ 
sufficiently big the natural ${\rm GL}(Q)$-equivariant pairing 
$\Gamma({\mathbb P}(Q),{\mathcal V}(m))\otimes{\mathcal O}(-m)
\longrightarrow{\mathcal V}$ is surjective. As ${\mathcal O}(-m)
\subset{\rm Sym}^mQ\otimes_k{\mathcal O}$, we can take the 
$k$-linear $G$-representation $\Gamma({\mathbb P}(Q),
{\mathcal V}(m))\otimes_k{\rm Sym}^mQ$ as $W$. \qed) 

\section{Admissible semi-linear representations} \label{main-conj} 
Let $G=G_{F/k}$ be the automorphism group of an algebraically 
closed field extension $F$ of $k$ of countable transcendence degree. 
This is a topological group with base of open subgroups $\{G_{F/k(x)}
~|~\mbox{for all $x\in F$}\}$, more details can be found in \cite{rep}. 

\subsection{A forgetful functor} \label{pgl-systems} 
Let $L_1\subset L_2\subset L_3\subset\dots$ be a sequence of purely 
transcendental extensions of finite type over $k$. Set $L_{\infty}
:=\bigcup_{j\ge 1}L_j$. Assume that $F$ is algebraic over $L_{\infty}$. 

Let ${\mathcal P}$ be the category with objects %, resp. ${\mathcal P}^+$,
$V_1\hookrightarrow V_2\hookrightarrow V_3\hookrightarrow\dots$, 
where \begin{itemize} \item $V_j$ is an 
$L_j\langle G_{L_j/k}\rangle$-%, resp. $L_j\langle{\rm End}(L_j/k)\rangle$-, 
module for each $j\ge 1$; 
\item $V_i\hookrightarrow V_j$ is a morphism of 
$L_i\langle G_{\{L_j,L_i\}/k}\rangle$-%, resp. $L_j\langle{\rm End}
%(\{L_j,L_i\}/k)\rangle$-, 
modules for all $i<j$ 
(here $G_{\{L_j,L_i\}/k}$ %, resp. ${\rm End}(\{L_j,L_i\}/k)$,
is the group of automorphisms %, resp. endomorphisms, 
of $L_j$ over $k$ inducing automorphisms %, resp. endomorphisms, 
of $L_i$); 
\item $V_s=V_i\cap V_j^{G_{L_j/L_s}}$ for any $s<i<j$. 
\end{itemize} 
The morphisms are defined as commutative diagrams 
$$\begin{array}{ccccccc} V_1 & \hookrightarrow & V_2 & 
\hookrightarrow & V_3 & \hookrightarrow & \dots \\ 
\varphi_1\downarrow\phantom{\varphi_1} & & 
\varphi_2\downarrow\phantom{\varphi_2} & & 
\varphi_3\downarrow\phantom{\varphi_3} \\ 
V_1' & \hookrightarrow & V_2' & \hookrightarrow & V_3' & 
\hookrightarrow & \dots \end{array}$$ where $\varphi_j$ 
%s an $L_j$-linear morphism of $\overline{H}_j$-modules. 
is an $L_j$-linear morphism of $G_{L_j/k}$-modules. 
Clearly, ${\mathcal P}$ is an additive $k$-linear category with kernels. 

\begin{lemma} \label{quasi-desc} Let $L_0\subset L_1\subset L_2\subset F$ 
be a pair of non-trivial purely transcendental extensions of finite 
transcendence degree. Let $H$ be a subgroup of $G_{\{F,L_2\}/L_0}$ 
projecting onto a subgroup of $G_{L_2/L_0}$ containing the permutation 
group of a transcendence basis of $L_2$ over $L_0$ extending a 
transcendence basis of $L_1$ over $L_0$. Then the subgroup $G'$ in $G$ 
generated by $G_{F/L_1}$ and $H$ coincides with $G_{F/L_0}$. \end{lemma}
{\it Proof.} If we could show that $G_{F/\xi(L_1)}\subset G'$ 
for any $\xi\in G_{F/L_0}$ then $G'$ would contain an open 
normal subgroup in $G_{F/L_0}$. It follows from the simplicity 
of $G_{F/\overline{L_0}}$ (Theorem 2.9 of \cite{rep}) that 
$G'=G_{F/L_0}$. Let $x_1,\dots,x_N$ be generators of $L_1$ over 
$L_0$, and $y_1,\dots,y_N$ be generators of $\xi(L_1)$ over $L_0$. 
By induction on $N$ we check that there exist $\sigma,\tau\in G'$ 
such that $\sigma x_1,\dots,\sigma x_N,\tau y_1,\dots,\tau y_N$ 
are algebraically independent over $L_2$. Then there is $\alpha
\in G_{F/L_2}$ such that $y_j=\tau^{-1}\alpha\sigma x_j$ for all 
$1\le j\le N$, and therefore, $G_{F/\xi(L_1)}=\tau^{-1}\alpha
\sigma G_{F/L_1}(\tau^{-1}\alpha\sigma)^{-1}\subset G'$. 

Let $N=1$. If $y_1\in\overline{L_1}$ then there is $\tau\in H$ 
such that $y'_1=\tau y_1\not\in\overline{L_1}$. 
If $y_1\not\in\overline{L_1}$ then there exists $\tau'\in G_{F/L_1}$ 
such that $\tau'y'_1\not\in\overline{L_2}$. 

Let $N>1$ and $y_1,\dots,y_{N-1}$ be algebraically independent 
over $L_2$. Suppose that the elements $y_1,\dots,y_N$ are 
algebraically dependent over $L_2$, i.e., $P(y_1,\dots,y_N)=0$ 
for some minimal polynomial 
$P\in\overline{L_2}[T_1,\dots,T_N]-\overline{L_0}[T_1,\dots,T_N]$ 
with a coefficient equal to 1. If $y_1,\dots,y_N$ are algebraically 
independent over $L_1$ then there is $\sigma\in G_{F/L_1}$ 
such that $\sigma y_1,\dots,\sigma y_N$ 
are algebraically independent over $L_2$. 

If $P\in\overline{L_1}[T_1,\dots,T_N]$ then there exists 
$\sigma\in H$ such that at least one of coefficients of the polynomial 
$\sigma P$ over $\overline{L_2}$ is outside of $\overline{L_1}$. 
As one of coefficients of $\sigma P$ is equal to 1, we see that 
$\sigma y_1,\dots,\sigma y_N$ 
are algebraically independent over $L_1$. \qed 

{\it Remark.} Similarly, one checks that $G_{L_{\infty}/L_i}$ and 
$G_{\{L_{\infty},L_j\}/L_s}$ generate $G_{L_{\infty}/L_s}$ for any $s<i<j$. 

For a semi-group $H$ of endomorphisms of an extension $K$ 
of $k$ denote by $\mathfrak{SL}(H;K)$ the (abelian) category 
of smooth semi-linear representations of $H$ over $K$. 
\begin{corollary} Sending $V$ to $V_{\infty}:=V^{G_{F/L_{\infty}}}$ 
defines a faithful functor $$\mathfrak{SL}(G;F)
\stackrel{\Phi}{\longrightarrow}
\mathfrak{SL}({\rm End}(L_{\infty}/k);L_{\infty}).$$ 
\end{corollary}
{\it Proof.} Let $H$ be the sub-semigroup of $G$ consisting of the 
elements inducing endomorphisms of $L_{\infty}$. As $H$ dominates 
${\rm End}(L_{\infty}/k)$, the latter semi-group acts on $V_{\infty}$. 

To check that $V_{\infty}$ is smooth, note that any $v\in V_{\infty}$ 
is fixed by the sub-semigroup $H$ generated by $G_{F/L_{\infty}}$ and 
${\rm End}(F/L)$ for some $L\subset F$ of finite type over $k$. There 
is $N\ge 1$ such that $L\subseteq\overline{L_N}$, so 
$H\supseteq\langle G_{F/L_{\infty}},{\rm End}(F/\overline{L_N})\rangle$. 
By Prop.2.14\footnote{Let $L_1$ and $L_2$ be subextensions of $k$ in 
$F$ such that $\overline{L_1}\bigcap\overline{L_2}$ is algebraic 
over $L_1\bigcap L_2$ and ${\rm tr.deg}(F/L_2)=\infty$. 
Then the subgroup in $G$ generated by $G_{F/L_1}$ and 
$G_{F/L_2}$ is dense in $G_{F/L_1\bigcap L_2}$.} 
of \cite{rep}, $H\supseteq{\rm End}(F/L_N)\supset
{\rm End}(\{F,L_{\infty}\}/L_N)$, i.e., the stabilizer of $v\in V_{\infty}$ 
contains the open sub-semigroup ${\rm End}(L_{\infty}/L_N)$. 

Clearly, $V_{\infty}\longmapsto V_{\infty}\otimes_{L_{\infty}}F$ and 
$V\longmapsto V^{G_{F/L_{\infty}}}$ define equivalences between 
the categories $\mathfrak{SL}({\rm End}(L_{\infty}/k);L_{\infty})$ and 
$\mathfrak{SL}(H;F)$ inverse to each other. \qed 

\begin{lemma} \label{P-equi-mo} 
The categories $\mathfrak{SL}(G_{L_{\infty}/k};L_{\infty})$ 
%resp. of ${\rm End}(L_{\infty}/k)$, ${\mathcal P}^+$, 
and ${\mathcal P}$ are equivalent. \end{lemma} 
{\it Proof.} By Remark after Lemma \ref{quasi-desc}, $V_{\infty}
\longmapsto(V_{\infty}^{G_{L_{\infty}/L_j}})_{j\ge 1}$ defines 
a functor to ${\mathcal P}$. 

To construct the inverse functor, one has to recover the 
$G_{L_{\infty}/k}$-action on the $L_{\infty}$-space 
$V_{\infty}:=\bigcup_{j\ge 1}V_j$, or equivalently, 
$G_{L_{\infty}/k}\times V_s\longrightarrow V_{\infty}$ for any 
$s\ge 1$. As $V_s\subset V_{s+1}\subset\dots$, we know the 
action of the elements of the closure of the set $T_s:=\bigcup_{j\ge s}
G_{\{L_{\infty},L_j\}/k}$ on $V_s$. 

For any $s\ge 1$, any $N\ge 1$ and any pair $x_1,\dots,x_N$ and 
$y_1,\dots,y_N$ of collections in $L_{\infty}$ complementable 
to sets of generators of $L_{\infty}/k$ there is $M\ge N$ and 
$\sigma\in G_{\{L_{\infty},L_M\}/k}$ such that $\sigma x_j=y_j$ 
for all $1\le j\le N$, and therefore, the closure of $T_s$ 
coincides with $G_{L_{\infty}/k}$. 

Any subfield of $L_{\infty}$ of finite type 
over $k$ is contained in some $L_j$, so 
$V_{\infty}=\bigcup_{j\ge 1}V_{\infty}^{G_{L_{\infty}/L_j}}$ 
for any smooth $G_{L_{\infty}/k}$-module $V_{\infty}$, and 
thus, our functors are mutually inverse. \qed

\begin{lemma} The composition of $\Phi$ with the equivalence of 
Lemma \ref{P-equi-mo} is given by $V\longmapsto(V_1\hookrightarrow 
V_2\hookrightarrow V_3\hookrightarrow\dots)$, where $V_j=V^{G_{F/L_j}}$. 
\end{lemma}
{\it Proof.} By Lemma 6.1 of \cite{rep}, $V^{G_{F/L_{\infty}}}= 
\bigcup_{L\subset L_{\infty}}V^{G_{F/L}}$, where $L$ runs over 
subfields of finite type over $k$. Then such $L$ are contained 
in appropriate $L_j$, so $V^{G_{F/L_{\infty}}}=
\bigcup_{j\ge 1}V^{G_{F/L_j}}$. 

The composition of $\Phi$ with the equivalence of Lemma \ref{P-equi-mo} 
is given by $V\mapsto(V_{\infty}^{G_{L_{\infty}/L_j}})_{j\ge 1}$. 
If $v\in V_{\infty}^{G_{L_{\infty}/L_s}}=V^{G_{\{F,L_{\infty}\}/L_s}}$ 
then there is $i\ge s$ such that $v\in V_i$. The subgroup 
$G_{\{F,L_{\infty},L_i\}/L_s}$ of $G_{\{F,L_{\infty}\}/L_s}$ 
acts on $V_i$ and its action factors through $G_{L_i/L_s}$, 
so $v\in V_i^{G_{L_i/L_s}}$. Then for any $j>i$ one has 
$v\in V_i^{G_{L_i/L_s}}\cap V_j^{G_{L_j/L_s}}$. 
Now Lemma \ref{quasi-desc} implies that $v\in V_s$, 
i.e., $V_{\infty}^{G_{L_{\infty}/L_s}}=V_s$.\qed 

\subsection{} \label{conj-cor} 
It is shown in \cite{rep} that the category of pure covariant 
motives over $k$ (defined using numerical equivalence) can be 
considered as a full subcategory of the category of graded 
semi-simple admissible $G$-modules of finite type. 
It is therefore desirable to know, 
whether the inverse inclusion holds. 

Any admissible representation $W$ of $G$ has the property 
that $W^{G_{F/L}}=W^{G_{F/L'}}$ for any extension $L$ of $k$ 
in $F$ and any purely transcendental extension $L'$ of $L$ 
in $F$. Such smooth representations form a Serre subcategory, 
denoted by ${\mathcal I}_G$, in the category of all smooth 
representations ${\mathcal S}m_G$ (cf. \S6 of \cite{rep}). 

The inclusion functor ${\mathcal I}_G\longrightarrow{\mathcal S}m_G$ 
admits the left adjoint ${\mathcal S}m_G\stackrel{{\mathcal I}}
{\longrightarrow}{\mathcal I}_G$. For an irreducible variety $Y$ 
over $k$ denote by $C_{k(Y)}$ the image under ${\mathcal I}$ of 
the free ${\mathbb Q}$-space with a basis given by the set of all 
embeddings of the function field $k(Y)$ into $F$ over $k$. 
Clearly, $C_{k(Y)}$ surjects onto $CH_0(Y_F)_{{\mathbb Q}}$ 
and there are several reasons to expect that $C_{k(Y)}=
CH_0(Y_F)_{{\mathbb Q}}$ when $Y$ is smooth and proper. 

\vspace{5mm} 

{\bf Conjecture.} 
\begin{enumerate} \item \label{irr-I-adm} Any irreducible object 
of ${\mathcal I}_G$ is contained in an admissible $F$-semi-linear 
representation $V$ of $G$, i.e., such that $\dim_{F^U}V^U<\infty$ 
for any open subgroup $U$ of $G$ (or equivalently, 
$\dim_LV^{G_{F/L}}<\infty$ for any subfield $L\subset F$ 
of finite type over $k$). 
\item \label{adm-Omega} Any irreducible admissible $F$-semi-linear 
representation of $G$ is contained in the tensor algebra 
$\bigotimes_F^{\bullet}\Omega^1_{F/k}$. \end{enumerate} 

{\it Remarks.} 1. To show the existence of smooth irreducible 
representations of $G$ not contained in the tensor algebra 
$\bigotimes^{\bullet}_F\Omega^1_{F/k}$, consider the $G$-module 
$W={\mathbb Q}[G/G_{\{F,\overline{k(x)}\}/k}]^{\circ}$, which 
is the same as the linear combinations of degree zero in the 
${\mathbb Q}$-space with the basis given by all algebraically 
closed subfields in $F$ of transcendence degree 1 over $k$. 

For any $x,y\in F-k$ algebraically independent over $k$, the 
vector $e=[\overline{k(x)}]-[\overline{k(y)}]$ is cyclic, since 
$\sigma e-e=[\overline{k(\sigma x)}]-[\overline{k(x)}]$ 
for any $\sigma\in G_{\{F,\overline{k(y)}\}/k}$. 
This means that $W$ admits an irreducible quotient. 
%For any $m\ge 1$, any $a_0,\dots,a_m\in{\mathbb Q}^{\times}$ with 
%$\sum_ja_j=0$ and any $x_0,\dots,x_m\in F-k$ such that the elements 
%of the collection $(x_1,\dots,x_m)$ and of the pairs $(x_0,x_j)$ 
%for all $1\le j\le m$ are algebraically independent over $k$, the 
%vector $e'=\sum_{j=0}^ma_j[\overline{k(x_j)}]$ is cyclic, since 
%$\sigma e'-e'=a_0([\overline{k(\sigma x_0)}]-[\overline{k(x_0)}])$ 
%for any $\sigma\in G_{\{F,\overline{k(x_1)},\dots,\overline{k(x_m)}\}/k}$. 
%This means that $W$ admits an irreducible quotient. 
%Let $e=[\overline{k(x)}]-[\overline{k(y)}]$ for some 
%$x,y\in F-k$ algebraically independent over $k$. 

Then ${\rm Hom}_G(W,W')\hookrightarrow(W')^{{\rm Stab}_e}$. 
As ${\rm Stab}_e=G_{\{F,\overline{k(x)},\overline{k(y)}\}/k}$, 
we see that $(W')^{{\rm Stab}_e}=
(\bigotimes^{\bullet}_{\overline{k(x,y)}}
\Omega^1_{\overline{k(x,y)}/k})
^{G_{\{\overline{k(x,y)},\overline{k(x)},\overline{k(y)}\}/k}}=k$, 
if $W'=\bigotimes^{\bullet}_F\Omega^1_{F/k}$, i.e., any morphism 
$W\longrightarrow W'$ factors through 
$W\stackrel{\varphi}{\longrightarrow}k\hookrightarrow W'$. 
Then $-\varphi(e)=\varphi(\sigma e)=\sigma\varphi(e)=\varphi(e)$ 
for any $\sigma\in G$ such that $\sigma x=y$ and $\sigma y=x$, 
and thus, $\varphi(e)=0$. This shows that no non-zero quotient of 
$W$ can be embedded into $\bigotimes^{\bullet}_F\Omega^1_{F/k}$. \qed 

2. It is known (cf. Prop.5.4 of \cite{rep}) that any smooth semi-linear 
representation of $G$ finite-dimensional over $F$ is trivial. 

%3. \item Any smooth semi-linear representation of $G$ carries a 
%canonical decreasing (weight) filtration with semi-simple quotients. 

%
%It seems that it is possible to deduce the first part of Corollary 
%\ref{adm-impl-adm} below from Conjecture (1), proceeding as follows. 
%
%Let $V$ be a smooth irreducible $F$-semi-linear
%$G$-representation, and $W\subset V$ be an irreducible 
%subrepresentation. Suppose that for an extension $L$ of finite 
%type over $k$ one has $\dim_{{\mathbb Q}}W^{G_{F/L}}=\infty$. 
%Let $Y$ be a smooth proper model of $L$ over $k$, and 
%${\mathcal V}_Y$ be the vector bundle on $Y$ corresponding to $V$. 
%If $\dim_k\Gamma(Y,{\mathcal V}_Y)<\dim_k{\rm Im}
%(W^{G_{F/L}}\otimes k\longrightarrow V)$ (in particular, if 
%$\dim_k{\rm Im}(W^{G_{F/L}}\otimes k\longrightarrow V)=\infty$ and 
%${\rm rk}{\mathcal V}_Y=\dim_LV^{G_{F/L}}<\infty$) then there is a 
%rational map $Y\stackrel{\pi}{\dasharrow}{\mathbb A}^{\dim Y}_k$ 
%defined in generic points of irreducible components of the divisor 
%of poles of some element $w\in W^{G_{F/L}}$ and separating these 
%components. Then $\pi_{\ast}$ does not annihilate $w$. On the other 
%hand, if $W$ is an object of ${\mathcal I}_G$ then the target of 
%$\pi_{\ast}$ is $W^G$, which is zero unless $W\cong{\mathbb Q}$. 

3. The following claim implies that the admissible semi-linear 
representations of $G$ form a tensor category. 
\begin{lemma} Let $E$ be either $F$ or any field of characteristic 
zero with the trivial $G$-action, and $W_1,W_2$ be smooth semi-linear 
representations of $G$ over $E$. Then, for any extension $L\subset F$ 
of $k$, one has 
$(W_1\otimes_EW_2)^{G_{F/\overline{L}}}=W_1^{G_{F/\overline{L}}}
\otimes_{E^{G_{F/\overline{L}}}}W_2^{G_{F/\overline{L}}}$. 

{\rm (This is not true if ${\rm tr.deg}(F/k)<\infty$. Namely, 
let $W_1$ and $W_2$ be of degree one, non-trivial and dual to 
each other. Then $(W_1\otimes_EW_2)^G=E^G$, but $W_1^G=W_2^G=0$.)} 
%
%Conversely, if for a subgroup $H\subset G$ one has 
%$(W_1\otimes_EW_2)^H=W_1^H\otimes_{E^H}W_2^H$ then 
%the field $F^H$ is algebraically closed. 
\end{lemma} 
{\it Proof.} Let $w=\sum_{j=1}^Nf_j\otimes g_j\in
(W_1\otimes_EW_2)^{G_{F/\overline{L}}}$ with minimal possible dimension 
of the $E$-subspace $H$ in $W_1$ spanned by $f_1,\dots,f_N$. 
Then $w\in H\otimes_EW_2$, and therefore, $w\in\sigma(H)\otimes_EW_2$ 
for any $\sigma\in G_{F/\overline{L}}$. Suppose that 
$f_j\not\in W_1^{G_{F/\overline{L}}}$ for some $j$ and that there 
is $\sigma\in G_{F/\overline{L}}$ such that $H\neq\sigma(H)$. 
Then $w\in(H\otimes_EW_2)\cap(\sigma(H)\otimes_EW_2)=H'\otimes_EW_2$, 
where $H'=H\cap\sigma(H)$ is a proper $E$-subspace in $H$, 
contradicting the minimality of dimension of $H$. This means that 
either $f_j\in W_1^{G_{F/\overline{L}}}$ for all $j$, i.e., 
$w\in W_1^{G_{F/\overline{L}}}\otimes_{E^{G_{F/\overline{L}}}}W_2$, 
or $H=\sigma(H)$ for any $\sigma\in G_{F/\overline{L}}$. In the first 
case one obviously has $w\in W_1^{G_{F/\overline{L}}}
\otimes_{E^{G_{F/\overline{L}}}}W_2^{G_{F/\overline{L}}}$. 

If $H=\sigma(H)$ for any $\sigma\in G_{F/\overline{L}}$ and 
${\rm tr.deg}(F/L)=\infty$ then $H$ is a smooth semi-linear representation 
of $G_{F/\overline{L}}$ of finite degree, which is trivial by Prop.5.4 of 
\cite{rep} in the case $E=F$ and by Theorem 2.9 of \cite{rep} if $E$ is a 
trivial $G$-module, i.e., $w\in H\otimes_EW_2\subset
(E\otimes_{E^{G_{F/\overline{L}}}}W_1^{G_{F/\overline{L}}})\otimes_EW_2$, 
and therefore, $w\in W_1^{G_{F/\overline{L}}}
\otimes_{E^{G_{F/\overline{L}}}}W_2^{G_{F/\overline{L}}}$. 

For arbitrary $L$ one has (cf. Lemma 6.1 of \cite{rep}) 
\begin{multline*} (W_1\otimes_EW_2)^{G_{F/\overline{L}}}=
\bigcup_{K\subset\overline{L}}(W_1\otimes_EW_2)^{G_{F/\overline{K}}}\\ 
=\bigcup_{K\subset\overline{L}}W_1^{G_{F/\overline{K}}}
\otimes_{E^{G_{F/\overline{K}}}}
W_2^{G_{F/\overline{K}}}=W_1^{G_{F/\overline{L}}}
\otimes_{E^{G_{F/\overline{L}}}}W_2^{G_{F/\overline{L}}},\end{multline*}
where $K$ runs over the set of subfields of finite type over $k$. 
%
%Suppose that for a subgroup $H\subset G$ one has 
%$(F\otimes F)^H=F^H\otimes F^H$. Set $L:=F^H$ then $H\subseteq G_{F/L}$ 
%and $L\otimes L=F^H\otimes F^H$ is a proper subspace in $(\overline{L}
%\otimes\overline{L})^{G_{\overline{L}/L}}=(F\otimes F)^{G_{F/L}}
%\subseteq(F\otimes F)^H$ the field $F^H$ is algebraically closed.
\qed   
%$W\longmapsto W^{G_{F/\overline{L}}}$ provides an equivalence 
%between the categories of smooth $G$- and $G_{\overline{L}/k}$-modules

\begin{proposition} \label{wedge} 
Let $W\in{\mathcal I}_G$ and $q\ge 0$ be an integer. 
Then \begin{itemize} \item any morphism $W\stackrel{\varphi}
{\longrightarrow}\left(\Omega^1_{F/k}\right)^{\otimes_F^q}$ 
factors through $W\longrightarrow\Omega^q_{F/k}\subseteq
\left(\Omega^1_{F/k}\right)^{\otimes_F^q}$; \item for any 
smooth proper variety $Y$ over $k$ an embedding $k(Y) 
\stackrel{\iota}{\hookrightarrow}F$ induces an injection 
$\varphi(W)\bigcap\iota_{\ast}\Omega^q_{k(Y)/k}\hookrightarrow
\Gamma(Y,\Omega^q_{Y/k})$ and there are the following 
canonical isomorphisms \begin{equation} \label{ch-0-forms} 
{\rm Hom}_G(C_{k(Y)},\bigotimes\nolimits^{\bullet}_F\Omega^1_{F/k})
\stackrel{\sim}{\longleftarrow}\Gamma(Y,\Omega^{\bullet}_{Y/k})
\stackrel{\sim}{\longrightarrow}{\rm Hom}_G(CH_0(Y_F),
\bigotimes\nolimits^{\bullet}_F\Omega^1_{F/k}).\end{equation} 
\end{itemize} \end{proposition}
{\it Proof.} Replacing $W$ with $\varphi(W)$, we may suppose that 
$W\subset\left(\Omega^1_{F/k}\right)^{\otimes_F^q}$. Let $\omega$ 
be a non-zero element of $W$. There is an extension $L\subset F$ 
of finite type over $k$ such that $\omega\in W^{G_{F/L}}$. Let 
$Y$ be a smooth projective model of $L$ over $k$. Then $\omega$ 
can be considered as a rational section of the coherent sheaf 
$\Omega^q_{Y^q/k}|_{\Delta_Y}$. 

If the divisor $(\omega)_{\infty}$ of poles of $\omega$ is 
non-zero, then there exists a generically finite rational map 
$f=(f_1,\dots,f_{\dim Y}):Y\dasharrow{\mathbb A}^{\dim Y}_k$ 
well-defined at generic points of irreducible 
components of $(\omega)_{\infty}$ and separating them. 

Then the direct image ${\rm tr}_{/k(f_1,\dots,f_{\dim Y})}
(\omega)$ of $\omega$ in $\Omega^q_{({\mathbb A}^{\dim Y}_k)^q/k}
|_{\Delta_{{\mathbb A}^{\dim Y}_k}}$ (in other words, in 
$\left(\Omega^1_{k(f_1,\dots,f_{\dim Y})/k}\right)^{
\otimes_{k(f_1,\dots,f_{\dim Y})}^q}$)
has poles, and in particular, it is non-zero. 

This means that ${\rm tr}_{/k(f_1,\dots,f_{\dim Y})}(\omega)$ is a 
non-zero element of $W^{G_{F/k(f_1,\dots,f_{\dim Y})}}=W^G$. 

On the other hand, $W^G\subseteq(\left(\Omega^1_{F/k}
\right)^{\otimes_F^q})^G=0$ if $q\ge 1$. 

Let
$\omega\in(\Omega^1_{L/k})^{\otimes_L^q}\backslash\Omega^q_{L/k}$ 
and $q\ge 2$. Fix a transcendence basis $x_1,\dots,x_{\dim Y}$ of 
the extension $L/k$. Then one has $\omega=
\sum_If_Idx_{i_1}\otimes\dots\otimes dx_{i_q}$. 

If there is a non-zero summand such that $dx_i$ appears $s\ge 2$ 
times then (after renumbering of the coordinates) we may suppose 
that $dx_1$ appears $s$ times. 

Otherwise, (after renumbering of the coordinates) we may suppose 
that one of the summands is $f_{1\dots q}dx_1\otimes\dots\otimes dx_q$ 
and $f_{123\dots q}\neq -f_{213\dots q}$. Replacing $x_2$ with 
$x_2+x_1$, we get the summand $-(f_{123\dots q}+f_{213\dots q})
dx_1\otimes dx_1\otimes dx_3\otimes\dots\otimes dx_q$. 

This means that $f_J\neq 0$, where $1$ appears in $J=(i_1\dots 
i_q)$ exactly $s\ge 2$ times. 

Let $L_0$ be the algebraic closure of $k(x_2,\dots,x_{\dim Y})$ 
in $L$, and $g$ be the genus of a smooth proper model $X$ of 
$L$ over $L_0$. Let $L_1\subset F$ be a finite extension of $L_0$ 
such that $X$ has a rational $L_1$-point $P$. Let $x_1$ be a local 
coordinate on $X_{L_1}$ in a neighbourhood of $P$ and $f_J(P)\neq 
0$. By Riemann--Roch theorem, for any effective reduced divisor 
$D$ on $X$ of degree $\ge 2g$ disjoint from $P$ there is a function 
$X_{L_1}\stackrel{z}{\longrightarrow}{\mathbb P}^1_{L_1}$ such that 
$z^{-1}(0)=s\cdot P+D$. 

Let $u$ be a local coordinate on ${\mathbb P}^1_{L_1}$ in the 
formal neighbourhood $V_P$ of $0$. There is a local coordinate $t$ 
on $X_{L_1}$ in a neighbourhood of $P$ such that $t^s=u\circ z$. 

Then $\omega$ is a sum of $g_Jdx_{i_1}\otimes\dots\otimes 
dx_{i_q}$, where $dt$ appears $s$ times and $g_J(P)\neq 0$, and 
some independent tensors. This implies that $\pi_{\ast}\omega$ 
is a sum of $(\pi|_{V_P})_{\ast}(g_J)s^{-s}u^{1-s}
dx_{i_1}\otimes\dots\otimes dx_{i_q}$, where $du$ appears $s$ 
times, some terms holomorphic at $0$, and some independent 
tensors. So $\pi_{\ast}\omega$ has a pole, and in particular, 
it is non-zero. Here $Z\stackrel{\pi}{\longrightarrow}
{\mathbb P}^1\times S$ is a thickening of $z$, $Z$ is a model 
of $LL_1$ over $k$, and $S$ is a model of $L_1$ over $k$. 

Both maps of (\ref{ch-0-forms}) are injective. Any morphism in 
question is determined by its value $\omega$ on the generator 
$\iota$ of $C_{k(Y)}$ (or of $CH_0(Y)$). As $\omega$ is fixed by 
$G_{F/\iota(k(Y))}$, the above implies that $\iota$ identifies 
$\omega$ with a regular differential form on $Y$. \qed 

\begin{lemma} For any integer $q\ge 0$ the semi-linear 
$G$-representation $\Omega^q_{F/k}$ is irreducible. \end{lemma}
{\it Proof.} Let $\omega$ be a non-zero element of $\Omega^q_{F/k}$. 
There is an extension $L$ of $k$ of finite type in $F$ such that 
$\omega\in\Omega^q_{L/k}$. Choose a transcendence basis 
$x_1,x_2,\dots$ of $L$ over $k$. 

Then $\omega=\sum_{I=(i_1<\dots<i_q)}f_Idx_I$, where $dx_I:=
dx_{i_1}\wedge\dots\wedge dx_{i_q}$. Let $\eta=\sum_{I=(i_1<\dots<i_q)}
g_Idx_I$ be a non-zero element in the semi-linear $G$-subrepresentation 
$\langle\omega\rangle_{F\langle G\rangle}$ generated by $\omega$ 
such that the number $N$ of non-zero $g_I$'s is minimal. 
If $N\ge 2$ and $g_J\neq 0$ then, for some $\sigma\in G$ with 
$\sigma x_i=\lambda_ix_i+\mu_i$, where $\lambda_i\in k^{\times}$, 
$\mu_i\in k$, $\prod_{i\in J}\lambda_i=1$, consider 
$\eta/g_J-\sigma(\eta/g_J)=\sum_{I\neq J}
\left(g_I/g_J-\sigma(g_I/g_J)\prod_{i\in I}\lambda_i\right)dx_I$. 

Then $\eta/g_J-\sigma(\eta/g_J)$ is a `shorter' non-zero element in 
$\langle\omega\rangle_{F\langle G\rangle}$. This contradicts our 
assumption, so $\langle\omega\rangle_{F\langle G\rangle}$ contains 
$dx_I$, and therefore, $\langle\omega\rangle_{F\langle G\rangle}
=\Omega^q_{F/k}$. \qed

\begin{corollary} \label{adm-impl-adm} If Conjecture holds 
then any irreducible object of ${\mathcal I}_G$ is admissible. 
If, moreover, numerical equivalence coincides with 
homological then there is an equivalence of categories 
$$\left\{\begin{array}{c} \mbox{{\rm covariant motives over $k$}}\\
\mbox{{\rm modulo numerical equivalence}}\end{array}\right\}
\stackrel{{\mathbb B}^{\bullet}}{\longrightarrow}
\left\{\begin{array}{c} \mbox{{\rm graded semi-simple admissible}}\\
\mbox{{\rm representations of} $G$ {\rm of finite length}}\end{array}
\right\}.$$ \end{corollary}
{\it Proof.} Let $W$ be an irreducible object of 
${\mathcal I}_G$. There is a smooth proper variety $Y$ 
over $k$ and a surjection $C_{k(Y)}\longrightarrow
\!\!\!\!\!\rightarrow W$. Assuming Conjecture, Proposition 
\ref{wedge} implies that $W$ can be embedded to $\Omega^q_{F/k}$ 
for an appropriate integer $q\ge 0$. As ${\rm Hom}_G(C_{k(Y)},
\Omega^q_{F/k})=\Gamma(Y,\Omega^q_{Y/k})$, any homomorphism 
$C_{k(Y)}\longrightarrow\Omega^q_{F/k}$ factors through 
$A^{\dim Y}(Y_F)$, where $A^{\ast}$ is the space of cycles 
``modulo (de Rham) homological equivalence over $k$''. More 
precisely, $A^{\ast}(Y_F)$ is the image of $CH_0(Y_F)_{{\mathbb Q}}$ 
in $H^{2\ast}_{dR/k}(Y_F)$. As $A^{\dim Y}(Y_F)$ is admissible, 
so is its quotient $W$. 

The fully faithful functor ${\mathbb B}^{\bullet}$ from 
\cite{rep} is given by the graded sum $\oplus_i
\lim\limits_{_L\longrightarrow}{\rm Hom}({\rm Prim}_L\otimes
{\mathbb L}^{\otimes i},-)$, where ${\mathbb L}$ is the Lefschetz 
motive, and ${\rm Prim}_L=\bigcap_{\varphi}\ker\varphi$ with 
$\varphi$ runing over all morphisms $Z\longrightarrow 
M\otimes{\mathbb L}$ for a smooth proper model $Z$ of $L$ over 
$k$ and effective motives $M$. 

It suffices to show that any irreducible admissible 
representation of $G$ is the degree-zero component of 
${\mathbb B}^{\bullet}(M)$ for a motive $M$. As $W$ is a 
quotient of $A^{\dim Y}(Y_F)$, this follows from the fact that 
$A^{\dim Y}(Y_F)$ coincides with the degree-zero component of 
${\mathbb B}^{\bullet}(Y)$, if numerical equivalence 
coincides with homological. \qed 

\vspace{5mm}

For a smooth proper $k$-variety $X$ denote by $C_{k(X)}^{\circ}$ 
the kernel of the morphism $C_{k(X)}\stackrel{\deg}{\longrightarrow}
{\mathbb Q}$, and by ${\mathcal F}^2C_{k(X)}$ 
the kernel of the morphism $C_{k(X)}^{\circ}\longrightarrow
{\rm Alb}X(F)_{{\mathbb Q}}$ induced by the Albanese map. 
\begin{corollary} \label{bl-con} If Conjecture holds then 
the folowing conditions are equivalent: \begin{itemize} \item 
the Albanese map identifies $CH_0(X_F)^0$ with ${\rm Alb}X(F)$; 
\item ${\mathcal F}^2C_{k(X)}=0$; \item $\Gamma(X,\Omega^q_{X/k})=0$ 
for all $q\ge 2$. \end{itemize} \end{corollary} 
{\it Proof.} By Corollary 6.24 of \cite{rep}, the cyclic 
$G$-module $C_{k(X)}$ is isomorphic to ${\mathbb Q}\oplus
{\rm Alb}X(F)_{{\mathbb Q}}\oplus{\mathcal F}^2C_{k(X)}$, so 
the $G$-module ${\mathcal F}^2C_{k(X)}$ is cyclic, and therefore, 
it admits an irreducible quotient $W$, which is non-zero if 
${\mathcal F}^2C_{k(X)}$ is non-zero itself. As $W\in{\mathcal I}_G$, 
Conjecture implies that there is an integer $q\ge 0$ 
and an embedding $W\hookrightarrow\Omega^q_{F/k}$. 

However, ${\rm Hom}_G(C_{k(X)},\Omega^q_{F/k})={\rm Hom}_G({\mathbb Q}
\oplus{\rm Alb}X(F)_{{\mathbb Q}},\Omega^q_{F/k})$ if $q\le 1$, so one 
has $q\ge 2$. This means that ${\rm Hom}_G(C_{k(X)},\Omega^q_{F/k})
=\Gamma(X,\Omega^q_{X/k})$ is non-zero for some $q\ge 2$. \qed 

\appendix
\section{Semi-linear representations of $PGL$ of degree one}
\label{Semi-PGL}
Fix an $n$-dimensional projective $k$-space ${\mathbb P}^n_k$ 
and a complement to a hyperlane ${\mathbb A}\subset{\mathbb P}^n_k$ 
with coordinates $x_1,\dots,x_n$. 
Let $G={\rm Aut}({\mathbb P}^n_k/k)\cong{\rm PGL}_{n+1}k$ be the group 
of automorphisms of ${\mathbb P}^n_k$ and $L=k({\mathbb P^n_k})$ 
be the function field of ${\mathbb P}^n_k$. 

Let $A=\left(A_{ij}\right)_{1\le i,j\le n+1}
\in{\rm PGL}_{n+1}k$ act on $L$ by $x_j\longmapsto
\frac{A_{1j}x_1+\dots+A_{nj}x_n+A_{n+1,j}}{A_{1,n+1}x_1+\dots+
A_{n,n+1}x_n+A_{n+1,n+1}}$. 

The aim of this section is to show that in characteristic zero any 
$L$-semi-linear $G$-representations of degree one is a `rational 
$L$-tensor power' of the space $\Omega^n_{L/k}$ of differential 
forms on $L$ over $k$ of top degree up to a character of 
a torsion quotient of $k^{\times}$. 

\begin{proposition} For any characteristic zero field $k$ 
one has $H^1(G,L^{\times}/k^{\times})={\mathbb Z}$. \end{proposition}
{\it Proof.} Let $k^n\cong U_0\subset G$ be the translation group 
of ${\mathbb A}$. First, by induction on $n$, we check that 
$H^1(U_0,L^{\times}/k^{\times})=0$, and in particular $H^1(U_0,k^{\times})
\stackrel{\sim}{\longrightarrow}H^1(U_0,L^{\times})$. 

Let $\lambda_0=(0,\dots,0,1)\in U_0$ and $U'_0=k^{n-1}\times\{0\}
\subset U_0$. For any collection $(f_{\lambda})_{\lambda}$ 
presenting a class in $H^1(U_0,L^{\times}/k^{\times})$ and any 
$\lambda,\mu\in U_0$ one has $f_{\lambda}(x+\mu)/f_{\lambda}(x)=
f_{\mu}(x+\lambda)/f_{\mu}(x)$. Multiplying $f_{\lambda}$ with 
rational functions of type $h(x+\lambda)/h(x)$ (which does not 
change the cohomology class), we may suppose that there are no 
pairs of irreducible components of the support of the divisor of 
$f_{\lambda_0}$ that differ by a translation by an integer 
multiple of $\lambda_0$. 

Then, for any $\mu\not\in\frac{1}{N}\lambda_0{\mathbb Z}$ with a 
sufficiently big integer $N$ there are no pairs of 
irreducible components of the support of the divisor of 
$f_{\lambda_0}(x+\mu)/f_{\lambda_0}(x)$ that differ by a translation 
by an integer multiple of $\lambda_0$, and therefore, 
$f_{\lambda_0}(x+\mu)/f_{\lambda_0}(x)=f_{\mu}(x+\lambda_0)/f_{\mu}(x)$ 
if and only if $f_{\lambda_0}(x+\mu)/f_{\lambda_0}(x)$ is constant, 
which means that $f_{\lambda_0}(x)$ is constant itself.

This implies also that for any $\mu\in U_0$ one has, 
$f_{\mu}(x)\in k(x_1,\dots,x_{n-1})^{\times}$, and thus, for 
any $\lambda\in\lambda_0\cdot k$ we get $f_{\lambda}(x+\mu)/
f_{\lambda}(x)=1$, i.e., $f_{\lambda}(x)\in k^{\times}$. 
By the induction assumption, there is some $g\in 
k(x_1,\dots,x_{n-1})^{\times}$ such that 
$f_{\lambda}(x)g(x)/g(x+\lambda)\in k^{\times}$ for all 
$\lambda\in U'_0$, and thus, 
$f_{\lambda}(x)g(x)/g(x+\lambda)\in k^{\times}$ for all
$\lambda\in U_0$, i.e., $H^1(U_0,L^{\times}/k^{\times})=0$. 

The Hochschild--Serre spectral sequence for the normal subgroup 
$U_0$ of the stabilizer $P_-$ of the hyperplane at infinity 
$\{(a_1:\dots:a_n:0)\}\subset{\mathbb P}^n_k$ in $G$ gives: 
$$E^{\ast,0}_2=H^{\ast}(P_-/U_0,H^0(U_0,L^{\times}/k^{\times}))=0\quad
\mbox{and}\quad E^{\ast,1}_2=H^{\ast}(P_-/U_0,H^1(U_0,L^{\times}/
k^{\times}))=0,$$ 
so we get $H^1(P_-,L^{\times}/k^{\times})=0$. 

As the stabilizer $P$ of the hyperplane $\{(0:a_1:\dots:a_n)\}
\subset{\mathbb P}^n_k$ in $G$ is conjugated to the subgroup $P_-$ of 
$G$, one has $H^1\left(P,L^{\times}/k^{\times}\right)=0$, so any 
element of $H^1(G,L^{\times}/k^{\times})$ can be presented 
by a cocycle sending any element of $P_-$ to $1$, and sending 
any element $A$ of $P$ to $f(x)^{-1}\cdot Af(x)$ for some 
$f(x)\in k(x_1,\dots,x_n)^{\times}/k^{\times}$. 

In the matrix form the subgroups $P_-$, $P$ and $P\cap P_-$ 
look, respectively,  as 
$$\left(\begin{array}{cccc} 1 & 0 & \dots & 0 \\ 
\ast & \ast & \dots & \ast \\ \dots & \dots & \dots & \dots \\ 
\ast & \ast & \dots & \ast \end{array}\right),\quad 
\left(\begin{array}{cccc} \ast & \dots & \ast & \ast \\ 
\dots & \dots & \dots & \dots \\ \ast & \dots & \ast & \ast \\ 
0 & \dots & 0 & 1 \end{array}\right),\quad 
\left(\begin{array}{ccccc} \ast & 0 & \dots & 0 & 0 \\ 
\ast & \ast & \dots & \ast & \ast \\ 
\dots & \dots & \dots & \dots & \dots \\ 
\ast & \ast & \dots & \ast & \ast \\ 0 & 0 & \dots & 0 & 1 
\end{array}\right).$$ 

In particular, 
$f(x)\in(L^{\times}/k^{\times})^{P\cap P_-}=x_1^{{\mathbb Z}}$, 
so $f(x)=x_1^m$ for some $m\in{\mathbb Z}$. 

This shows that $f_A=(A_{1,n+1}x_1+\dots
+A_{n,n+1}x_n+A_{n+1,n+1})^{-m_1}$ for any $A\in P$. 

Since the subgroups $P_-$ and $P$ generate ${\rm PGL}_{n+1}k$, one has 
$$H^1(G,L^{\times}/k^{\times})=\left\{f_A=
(A_{1,n+1}x_1+\cdots+A_{n,n+1}x_n+A_{n+1,n+1})^{-m}
~|~m\in{\mathbb Z}\right\}\cong{\mathbb Z}.\quad\qed $$ 

\begin{corollary} \label{pgl-2} For any field $k$ of characteristic 
zero there is a natural exact sequence \begin{multline*} 
0\longrightarrow{\rm Hom}(k^{\times}/(k^{\times})^{n+1},k^{\times})
\longrightarrow 
H^1(G,L^{\times})\stackrel{p}{\longrightarrow}H^1(G,L^{\times}/
k^{\times})={\mathbb Z} \\ \longrightarrow{\rm image}\left[
{\mathbb Z}\longrightarrow{\rm End}(k^{\times})\otimes
{\mathbb Z}/(n+1){\mathbb Z}\right]\longrightarrow 0.\end{multline*} 

The class\footnote{of the cocycle $(\sigma\longmapsto\sigma
\omega/\omega)\in H^1(G,L^{\times})$ for any non-zero $n$-form 
$\omega\in\Omega^n_{L/k}$} of $\Omega^n_{L/k}$ generates 
a subgroup of index $n+1$ in $H^1(G,L^{\times}/k^{\times})$. 

In particular, the number of $(n+1)$-st roots of unity 
in $k$ divides the order of ${\rm coker}p$ (thus, $\Omega^n_{L/k}$ 
generates $H^1(G,L^{\times})$ if $k$ contains all $(n+1)$-st roots 
of unity and $k^{\times}$ is $(n+1)$-divisible). \end{corollary}

{\sc Example.} If $k={\mathbb R}$ then $p$ is bijective 
for even $n$, and $\#{\rm coker}p=\#\ker p=2$ for odd $n$. 

{\it Proof.} Since the commutant of $G$ coincides with 
${\rm PSL}_{n+1}k$, the determinant induces an isomorphism 
${\rm Hom}(k^{\times}/(k^{\times})^{n+1},k^{\times})
\stackrel{\sim}{\longrightarrow}{\rm Hom}(G,k^{\times})$, 
so the short exact sequence 
$$1\longrightarrow k^{\times}\longrightarrow L^{\times}
\longrightarrow L^{\times}/k^{\times}\longrightarrow 1$$ 
gives the map $p$ and determines its kernel. To identify its 
cokernel, suppose that the 1-cocycle 
$$A=\left(A_{ij}\right)_{1\le i,j\le n+1}
\longmapsto(A_{1,n+1}x_1+\dots+A_{n,n+1}x_n+A_{n+1,n+1})^{-m}$$ 
on $G$ with values in $L^{\times}/k^{\times}$ can be lifted to 
a 1-cocycle $A=\left(A_{ij}\right)_{1\le i,j\le n+1}\longmapsto
\Phi(A)\cdot(A_{1,n+1}x_1+\dots+A_{n,n+1}x_n+A_{n+1,n+1})^{-m}$ on 
$G$ (considered as a 1-cocycle on $GL(Q)$ for an $(n+1)$-dimensional 
$k$-vector space $Q$) with values in $L^{\times}$. 
Then $\Phi:GL(Q)\longrightarrow k^{\times}$ is a homomorphism, and 
thus, $\Phi$ factors through the determinant: $\Phi(A)=\phi(\det A)$ 
for a homomorphism 
$k^{\times}\stackrel{\phi}{\longrightarrow}k^{\times}$. The cocycle 
on $GL(Q)$ defined by $\Phi$ descends to a cocycle on $G$ if and 
only if $\Phi$ is homogeneous of degree $m$, so 
$\phi(\lambda)^{n+1}=\lambda^m$. This implies that $m$, considered 
as element of ${\rm End}(k^{\times})\supset{\mathbb Z}$, 
should be divisible by $n+1$. 

As any endomorphism of $k^{\times}$ induces an endomorphism of the 
subgroup of $(n+1)$-st roots of unity, if $k$ contains $t$ out of 
$n+1$ roots of unity of order $n+1$, then $n+1$ divides $m$ as 
element of ${\mathbb Z}/t{\mathbb Z}$, which simply means that 
$m\equiv 0\pmod t$. \qed

\section{Semi-linear representations of 
degree one of a subgroup of the Cremona group}
As in the previous section, we fix an $n$-dimensional projective 
$k$-space ${\mathbb P}^n_k$ and some affine coordinates 
$x_1,\dots,x_n$ on ${\mathbb P}^n_k$. 

Let $P={\rm Aut}({\mathbb P}^n_k/k)\cong{\rm PGL}_{n+1}k$ be the 
automorphism group of ${\mathbb P}^n_k$ (denoted by $G$ in Appendix 
\ref{Semi-PGL}), and $L=k({\mathbb P}^n_k)$ be the function field 
of ${\mathbb P}^n_k$. 

Let $G$ be the subgroup of the Cremona group ${\rm Cr}_n(k)=
{\rm Aut}(L/k)$ generated by $P$ and by the involution $\sigma$ such 
that $\sigma x_1=x_1^{-1}$ and $\sigma x_j=x_j$ for all $2\le j\le n$. 

The aim of this section is to show that in characteristic zero any 
$L$-semi-linear $G$-representations of degree one is an integral 
$L$-tensor power of the space $\Omega^n_{L/k}$ of differential 
forms on $L$. 

\begin{proposition} \label{cre} 
Let $k$ be a field of characteristic zero. Suppose that 
either $n$ is even, or $k^{\times}$ is 2-divisible. 
Then the isomorphism class of $\Omega^n_{L/k}$ generates 
the group $H^1(G,L^{\times})$. \end{proposition} 
{\it Proof.} Let $(a_{\tau})$ be a 1-cocycle on $G$. We may suppose 
that (in notation of the proof of Corollary \ref{pgl-2}) the 
restriction of $(a_{\tau})$ to $P$ coincides with $A\longmapsto\phi
(\det A)\cdot(A_{1,n+1}x_1+\dots+A_{n,n+1}x_n+A_{n+1,n+1})^{-m}$ for 
a homomorphism $k^{\times}\stackrel{\phi}{\longrightarrow}k^{\times}$. 

Let $T\subset P$ be the maximal torus subgroup such that 
$\tau x_j/x_j=:\lambda_j(\tau)\in k^{\times}$ for any 
$\tau\in T$ and any $1\le j\le n$. Then $\sigma$ normalizes 
$T$ and for any $\tau\in T$ one has $\lambda_j(\sigma^{-1}\tau\sigma)
=\lambda_j(\tau)^{1-2\delta_{1,j}}$, and 
$a_{\tau}=\phi(\lambda_1(\tau)\cdots\lambda_n(\tau))$. 
As $a_{\sigma^{-1}\tau\sigma}=
\phi(\lambda_1(\tau)^{-1}\cdot\lambda_2(\tau)\cdots\lambda_n(\tau))$ 
and $a_{\tau}=a_{\sigma}\cdot\sigma a_{\sigma^{-1}\tau\sigma}\cdot
\tau a_{\sigma}^{-1}$, this implies that $\tau a_{\sigma}^{-1}\cdot 
a_{\sigma}=\phi(\lambda_1(\tau)^2)$. 

This means, in particular, that $a_{\sigma}$ does not depend on the 
variables $x_2,\dots,x_n$, i.e., $a_{\sigma}\in k(x_1)^{\times}$, 
and $a_{\sigma}(\lambda_1x_1)=\phi(\lambda_1^{-2})\cdot 
a_{\sigma}(x_1)$ for any $\lambda_1\in k^{\times}$. 

It is now clear that $a_{\sigma}(x_1)$ is homogeneous, say of degree 
$s\in{\mathbb Z}$, so $\phi(\lambda^2)=\lambda^{-s}$. Evaluating both 
sides at $-1$, we see that $s$ is even. Recall from the proof of 
Corollary \ref{pgl-2} that $\phi(\lambda)^{n+1}=\lambda^m$, so 
$\lambda^{-s(n+1)}=\phi(\lambda^2)^{n+1}=\lambda^{2m}$, and thus, 
$m=-\frac{s}{2}(n+1)$ is divisible by $n+1$. Then 
%and $\phi(\lambda_1)=\lambda_1^{m/(n+1)}$ 
$a_{\sigma}(x_1)=c\cdot x_1^{-2m/(n+1)}$ for some $c\in k^{\times}$, 
%As $a_{\sigma}\cdot\sigma a_{\sigma}=1$, we get $c=\pm 1$. 
so $(a_{\tau})$ is the product on an integer power of the class 
of $\Omega^n_{L/k}$ and a homomorphism 
$G\stackrel{c}{\longrightarrow}k^{\times}$. 

One has $c:P\stackrel{\det}{\longrightarrow}k^{\times}/
(k^{\times})^{n+1}\stackrel{\phi}{\longrightarrow}k^{\times}$. As 
$\sigma\tau\sigma=\tau^{-1}$ for any $\tau={\rm diag}(\lambda,1,
\dots,1)\in P$, we get $c(\tau)^2=\phi(\lambda)^2=1$ for any 
$\lambda\in k^{\times}$. If either $n$ is even, or $k^{\times}$ 
is 2-divisible, this implies that $\phi\equiv 1$, so $c(P)=\{1\}$. 

Let $\iota_{01}$ be the involution in $P$ sending $(x_1,\dots,x_n)$ to 
$(1/x_1,x_2/x_1,\dots,x_n/x_1)$, $s_1=\iota_{01}\sigma\iota_{01}:
(x_1,\dots,x_n)\longmapsto(1/x_1,x_2/x_1^2,\dots,x_n/x_1^2)$ 
and $s_0:(x_1,\dots,x_n)\longmapsto(x^{-1}_1,\dots,x_n^{-1})$. 
The element $s_0$ belongs to $G$, since it is the product of the 
elements $\iota_{1j}\sigma\iota_{1j}$ for all $1\le j\le n$, where 
$\iota_{1j}$ are involutions in $P$ such that $\iota_{1j}x_s=x_s$ 
for $s\not\in\{1,j\}$ and $\iota_{1j}x_1=x_j$. 
Let $g_0$ be the involution in $P$ sending $(x_1,\dots,x_n)$ to 
$(\frac{x_1}{x_1-1},\frac{x_2}{x_1-1},\dots,\frac{x_n}{x_1-1})$. 
Then one has the following well-known identity in $G$: 
$s_1=g_0s_0g_0s_0g_0$. Then for the homomorphism $c$ as above 
one has $c(\sigma)=c(s_1)=c(g_0)^3c(s_0)^2=c(s_0)^2$. 
As $s_0^2=1$, this implies that $c(\sigma)=1$. \qed 

\vspace{5mm}

{\it Remarks.} 1. Let $\tilde{G}$ be the subgroup of 
${\rm Cr}_n(k)$ generated by $G$ and the involution $\xi$ such 
that $\xi x_1=x_1^{-1}$, $\xi x_2=x_1^{-1}x_2$ and $\xi x_j=x_j$ 
for any $3\le j\le n$. Then Proposition \ref{cre} remains valid if 
$G$ is replaced by $\tilde{G}$. Almost the same proof goes through. 

2. If $k$ is algebraically closed and $n=2$
then by M.Noether theorem $G={\rm Cr}_2(k)$.

\begin{proposition} \label{no-rep} Let $k$ be an algebraically 
closed field of characteristic zero, ${\mathcal A}$ a 
noetherian algebraic group scheme over a ring $R$ and $n\ge 2$. 
Then ${\rm Hom}(\tilde{G},{\mathcal A}(R))=\{1\}$. \end{proposition} 
{\it Proof.} It was shown at the end of the proof of Proposition 
\ref{cre} that there are no proper normal subgroups of $G$ containing 
$P$. As $P$ is simple (generated by any non-trivial conjugacy 
class), there are no proper normal subgroups of $G$ containing 
a non-trivial element of $P$. 

Let the elements of ${\rm Cr}_2(k)$ act identically on 
$k(x_3,\dots,x_n)$. Then ${\rm Cr}_2(k)$ becomes a subgroup 
in ${\rm Cr}_n(k)$. By M.Noether theorem, ${\rm Cr}_2(k)$ 
is generated by $\sigma$, $\xi$ and $A_2:=P\bigcap{\rm Cr}_2(k)$. 

Denote by $H\cong k(x_2)\rtimes k^{\times}$ 
the subgroup of ${\rm Cr}_2(k)\subseteq\tilde{G}$ consisting of elements 
$\tau=(q(x_2),b)$ such that $\tau x_1=x_1+q(x_2)$ and $\tau x_2=b\cdot x_2$ 
for some $q(x_2)\in k(x_2)$ and $b\in k^{\times}$. 

For any $N\ge 3$ and any primitive $N!$-th root of unity 
$\zeta_{N!}$ the centralizer of $(0,\zeta_{N!})$ in $H$ 
is $k(x_2^{N!})\rtimes k^{\times}$.

Let $\rho:\tilde{G}\longrightarrow{\mathcal A}(R)$ be a homomorphism. 
First, we show that $\ker\rho\bigcap H\bigcap A_2\neq\{1\}$. 

Suppose that $\ker\rho\bigcap H=\{1\}$. Then $H
\stackrel{\rho}{\hookrightarrow}{\mathcal A}(R)$, and thus, the 
centralizer of $(0,\zeta_{N!})$ in $H$ is the intersection of $H$ with 
the centralizer of $(0,\zeta_{N!})$ in ${\mathcal A}(R)$. The centralizer 
of an element of ${\mathcal A}(R)$ is the group of $R$-points of a closed 
subgroup in ${\mathcal A}$, so any descending sequence of centralizers 
should stabilize. This is not the case for the sequence 
$(k(x_2^{N!})\rtimes k^{\times})_{N\ge 1}$. 

Let $(q_1(x_2),b)\in\ker\rho\bigcap H-\{1\}$. If $b\neq 1$ then 
\begin{multline*} (x_2,1)(q_1(x_2),b)(x_2,1)^{-1}(q_1(x_2),b)^{-1}=
(x_2+q_1(x_2),b)(-x_2,1)(-q_1(b^{-1}x_2),b^{-1}) \\ 
=((1-b)x_2+q_1(x_2),b)(-q_1(b^{-1}x_2),b^{-1})=((1-b)x_2,1)
\in\ker\rho\bigcap H-\{1\},\end{multline*} 
so there is $(q(x_2),1)\in\ker\rho\bigcap H-\{1\}$. 

It is easy to see using prime decomposition of $q$ that 
$(x_1,x_2)\stackrel{\alpha}{\longmapsto}(q(x_2)x_1,x_2)$ 
is an element of $\tilde{G}$ (or, by M.Noether theorem, 
$\alpha\in{\rm Cr}_2(k)\subseteq\tilde{G}$), so 
$\alpha(q(x_2),1)\alpha^{-1}\in\ker\rho\bigcap H$. 

But $(x_1,x_2)\stackrel{\alpha(q,1)\alpha^{-1}}{\longmapsto}(x_1+1,x_2)$ 
is a non-trivial element of $A_2$, so $\ker\rho\bigcap H\bigcap 
A_2\neq\{1\}$. 

As $\ker\rho$ is a normal subgroup in $\tilde{G}$, this implies that
$\ker\rho\supseteq G$, so the image of $\rho$ is $\{1,\rho(\xi)\}$. 
The Euclid algorithm shows that $\sigma$, $\xi$ and the permutations
of the set $\{x_1,\dots,x_n\}$ generate a subgroup isomorphic to 
${\rm GL}_n{\mathbb Z}$. The commutators and $\sigma$ generate 
a subgroup contaning the congruence subgroup $\ker({\rm GL}_n{\mathbb Z}
\longrightarrow{\rm GL}_n{\mathbb F}_2)$, where ${\mathbb F}_2$ is the
field of 2 elements, so the image of $\rho$ is isomorphic to the image
of some homomorphism ${\rm SL}_n{\mathbb F}_2\longrightarrow\{\pm 1\}$,
sending $\sigma$ to 1. For any $n\ge 3$ the group ${\rm SL}_n
{\mathbb F}_2$ is simple, so $\ker\rho=\tilde{G}$, \qed

\vspace{5mm} 

\noindent
{\sl Acknowledgement.} {\small A part of this work was done at the 
Institut de Math{\'e}matiques de Luminy during the Semester 
`Arithm{\'e}tique et Th{\'e}orie de l'Information' in Spring 2003. 
I am grateful to its organizers Yves Aubry, Gilles Lachaud, 
Mikhail Tsfasman for the invitation and to CNRS for financial 
support in April--June 2003. 
%I would like to thank also Alexei Bondal for several enightening 
%discussions and for bringing \cite{K} to my attention. 

}

\vspace{5mm}

$$\begin{array}{l} \mbox{Independent University of Moscow} \\
\mbox{121002 Moscow} \\ \mbox{B.Vlasievsky Per. 11} 
\\ \mbox{{\tt marat@mccme.ru}} \end{array}\quad\mbox{and}\quad
\begin{array}{l}\mbox{Institute for Information} \\ 
\mbox{Transmission Problems} \\ 
\mbox{of Russian Academy of Sciences} \end{array}$$
\end{document}